\documentclass[12pt,leqno,oneside]{amsart}

\usepackage{tikz-cd}
\usepackage[utf8]{inputenc}

\definecolor{Violet}{RGB}{120,40,200}
\newcommand{\Vio}{
\textcolor{Violet}}

\usepackage{enumitem}
\usepackage{array}
\usepackage{amsmath,amstext,amsthm,amssymb,amscd}
\usepackage{amsfonts}

\usepackage{mathrsfs,dsfont}
\usepackage{typearea}
\usepackage{pdfsync}
\usepackage{mathtools}
\usepackage
{hyperref}
\numberwithin{equation}{section}

\newtheorem{Theorem}{Theorem}[section]
\newtheorem{Proposition}[Theorem]{Proposition}
\newtheorem{Lemma}[Theorem]{Lemma}
\newtheorem{Corollary}[Theorem]{Corollary}
\theoremstyle{definition}

\newtheorem{Remark}[Theorem]{Remark}

\usepackage[math]{anttor}
\usepackage{pifont}

\usepackage{MnSymbol,wasysym}

\newcommand{\bs}{\backslash}
\newcommand{\R}{\mathbb{R}}

\newcommand{\Z}{\mathbb{Z}}

\newcommand{\eproof}{
\begin{flushright}
\Vio{\ding{170}}
\end{flushright}}

\newcommand{\ii}{\hspace*{1cm}}

\usepackage{setspace}
\usepackage{ragged2e}
\usepackage{graphicx}
\usepackage{geometry}
\usepackage{float}
\DeclareMathAlphabet{\mathdutchcal}{U}{dutchcal}{m}{n}

\begin{document}
\raggedright
\title[\resizebox{5.5in}{!}{Asymptotics of Fubini-Study Currents for Sequences of Line Bundles}]{{Asymptotics of Fubini-Study Currents for Sequences of Line Bundles}}
\author{Melody Wolff}
\address{Department of Mathematics, Syracuse University, Syracuse, NY 13244-1150, USA}
\email{mawolff@syr.edu}

\subjclass[2010]{Primary 32L10; Secondary 31C10, 32A36, 32L05, 32Q15, 32U05, 32U40}
\keywords{Bergman kernel function, Fubini-Study current, 
singular Hermitian metric, big line bundle}

\date{October 6, 2024}

\pagestyle{myheadings}

\begin{abstract}
    We study the Fubini-Study currents and equilibrium metrics of continuous Hermitian metrics on sequences of holomorphic line bundles over a fixed compact K{\"a}hler manifold. We show that the difference between the Fubini-Study currents and the curvature of the equilibrium metric, when appropriately scaled, converges to $0$ in the sense of currents. As a consequence, we obtain sufficient conditions for the scaled Fubini-Study currents to converge weakly.
\end{abstract}

\maketitle
\tableofcontents

\section{Introduction}\label{S1} 

\ii In this paper, we will be working with sequences of holomorphic line bundles $\{L_p\}$ with continuous Hermitian metrics $h_p$, $p\geq 1$, and will be studying asymptotic properties of the Fubini-Study current first explored in \cite{CMM17}. We restrict our work to compact K{\"a}hler manifolds while allowing the metrics to have non-positive curvature. Although we won't require positivity conditions on our metrics, we will require the existence of positively curved metrics with growth conditions similar to those used in \cite{CMM17}.

\ii In 1988 Tian explored the case where $(L_p,h_p)=(L^p,h^p)$, with $(L^p,h^p)=(L^{\otimes p},h^{\otimes p})$ for some holomorphic line bundle $L$ equipped with a smooth metric $h$ (see \cite{T88}). He showed that if $(X,\omega)$ is a compact K{\"a}hler manifold with a line bundle $(L,h)$ such that the curvature $c_1(L,h)$ is positive and $h$ is smooth, then the normalized Fubini-Study forms $\gamma_p/p$ (see \ref{S2.3} for definition) converge to $c_1(L,h)$ in $C^2$. Later results showed that the convergence was actually in $C^\infty$  \cite{Ca99,Ru60,Zel98}. We refer to \cite{MM07}. As a consequence of this result, Tian showed that particular K{\"a}hler-Einstein forms could be approximated by Fubini-Study currents, which answered a question of Yau \cite{Y86}. 

\ii Also included in \cite{MM07} is an asymptotic expansion of the Bergman kernel (see \cite[Theorem 5.4.10]{MM07}). This expansion has been shown to provide information about the underlying Kähler manifold. In particular, the asymptotic expansion can be used to prove the Kodaira embedding theorem (see \cite[section 5.1.2]{MM07}).

\ii The assumptions of Tian's results on $(L,h)$ were relaxed by Coman and Marinescu in \cite{CM15}. They worked in the case $(L_p,h_p)=(L^p,h^p)$ and showed that if $c_1(L,h)$ was an integrable K{\"a}hler current, then  the aforementioned convergence result holds in the sense of currents. The results were further generalized by Coman, Ma, and Marinescu (see \cite{CMM17}). They showed that if $c_1(L_p,h_p)\geq a_p\omega$ where $a_p\to \infty$, and
\smallskip
\[A_p=\int_X c_1(L_p,h_p)\wedge\omega^{n-1},\]

then 
\[\frac{\gamma_p-c_1(L_p,h_p)}{A_p}\to 0\] 

weakly as currents.

\ii Berman's work in 2009 introduced the notion of an equilibrium metric $h^{eq}$ (See section \ref{S2.2} for definition) corresponding to a smooth metric $h$ on a holomorphic line bundle $L$ (see \cite{Ber09}). He worked in the setting where $(L_p,h_p)=(L^p,h^p)$. He showed that for any compact subset $\Omega$ of $X\bs \mathbb{B}_+(L)$, where $\mathbb{B}_+(L)$ is the augmented base locus (definition given in \cite{Ber09}), there exists $C_\Omega\geq 0$ such that 
\smallskip
\[-\frac{C_\Omega}{p}\leq \frac{\log P_p}{p}-(\varphi^{eq}-\varphi)\leq \frac{C_\Omega+n\log p}{p},\]

where $\varphi$ and $\varphi^{eq}$ are the global weights of $h$ and $h^{eq},$ respectively (see \cite[Theorem 1.5]{Ber09}). Another interesting result of his, is the convergence

\[\limsup_p \frac{\dim H^0(X,L^p)}{p^n}=\int_{U(L)} \frac{c_1(L,h^{eq})^n}{n!},\]

where $n=\dim(X)$ and $U(L)$ is the set where the weights of $h^{eq}$ are locally bounded. He also showed the following weak convergence of measures

\[\frac{P_p \omega^n}{p^n}\to \chi_{U(L)}\left(\frac{c_1(L,h^{eq})^n}{n!}\right),\]

where $P_p$ is the Bergman kernel function (as defined in Section \ref{S2.3}) and $\chi_{U(L)}$ is the characteristic function. 

\ii In 2019, Coman, Marinescu, and Nguy{\^e}n used the equilibrium metric to generalize Tian's work (see \cite[Cor. 5.7]{CMN19}). Like Tian and Berman, they worked in the case where $(L_p,h_p)=(L^p,h^p)$ and showed
\smallskip
\[
\frac{\gamma_p}{p}\to c_1(L,h^{eq})
\]

weakly as currents. 

\ii The following conditions will serve as the setting for most of our results in this paper:

\ii (A) $(X,\omega)$ is a compact (connected) K{\"a}hler manifold of complex dimension $n$.\\
\vspace{2mm}
\ii (B) $L_p$, $p\geq 1$, is a holomorphic line bundle on $X$ equipped with a continuous metric $h_p$ and a singular metric $g_p$ verifying
\smallskip
\begin{equation}\label{Defn ap}
    c_1(L_p,g_p)\geq a_p\omega\text{ on }X,\text{ where }\lim_{p\to\infty} a_p=\infty.
\end{equation}

Set 
\[A_p=\int_X c_1(L_p,g_p)\wedge \omega^{n-1}.\]\\
\medskip
\ii (B') $L_p$, $p\geq 1$, is a pseudo-effective holomorphic line bundle equipped with the continuous metric $h_p$. There exists an open coordinate polydisc cover $\{V_\alpha\}$ of $X$, frames $e_p^\alpha$ of $L_p$ on $V_\alpha$, functions $\phi^\alpha\in C(V_\alpha)$, and constants $A_p>0$, such that $A_p\to\infty$ and 
\[\phi_p^\alpha/A_p\to \phi^\alpha\text{ locally uniformly,}\]
where $\phi_p^\alpha$ is the local weight of $h_p$ corresponding to $e_p^\alpha.$\\
\vspace{2mm}
\ii Following Berman, we define $h_p^{eq}$ to be the equilibrium metric of $h_p$. As well, we recall that $\gamma_p$ is the Fubini-Study current. For definitions refer to sections \ref{S2.2} and \ref{S2.3}.

\ii For any coordinate polydisc $U$ and $p\geq 1$, let $e_p$ be a local frame of $L_p$ on $U$ (see Section \ref{S2.2}). Let $\phi_p:U\to\R$ be the continuous function such that
\[h_p(e_p,e_p)= e^{-2\phi_p}.\]

Similarly, we define $\rho_p:U\to[-\infty,\infty)$ by
\smallskip
\[g_p(e_p,e_p)= e^{-2\rho_p}.\]

The functions $\phi_p$ and $\rho_p$ are called the local weights of $h_p$ and $g_p$ on $U$.

    \ii Our main results are:

\begin{Theorem}\label{T1.1}
    Let $(X,\omega)$ and $(L_p,h_p)$ be as in $(A)$ and $(B)$. If every $x\in X$ has a neighborhood $U$ with local frames $e_p$ of $L_p$, such that the families of local weights $\{\phi_p/A_p\}$ and $\{\rho_p/A_p\}$ are uniformly bounded in $L^1(U)$, and there exists $M>0$ such that 
    \begin{equation}\label{T1.1-Inequality}
    -M\omega\leq  c_1(L_p,h_p)/A_p\leq M\omega,
    \end{equation}
    then
    \begin{equation}\label{Main-Convergence}
    \frac{\gamma_{p}-c_1(L_p,h_p^{eq})}{
    A_p}\to 0 \text{ weakly as currents.}
    \end{equation}
    
\end{Theorem}

\ii If $h_p$ verifies (\ref{Defn ap}), then $h_p=h_p^{eq}$. In this case, (\ref{Main-Convergence}) becomes

\[
\frac{\gamma_p-c_1(L_p,h_p)}{A_p}\to 0
\]

weakly as currents. This is a special case of the convergence condition shown in \cite{CMM17}. If in addition we assume that $(L_p,h_p)=(L^p,h^p)
$ and set $A_p=p$, our assumptions are automatically satisfied, and we obtain the convergence shown in \cite{CM15}.

\begin{Theorem}\label{T1.2}

Let $(X,\omega)$ and $(L_p,h_p)$ be as in $(A)$ and $(B)$. If every $x\in X$ has a neighborhood $U$ with local frames $e_p$, such that the collection of scaled local weights $\{\phi_p/A_p\}$ is equicontinuous and uniformly bounded, and $\{\rho_p/A_p\}$ is uniformly bounded in $L^1(U)$, then (\ref{Main-Convergence}) holds.

\end{Theorem}

\ii When $(L_p,h_p)=(L^p,h^p)$ equicontinuity is trivial, as in that case we can take $A_p=p$, and for particular local frames we have $\phi_p/p=\phi_1$. Like Theorem \ref{T1.1}, this result can be considered a generalization of the convergence in \cite{CM15}. Another case where equicontinuity is automatically satisfied is when $(L_p,h_p)$ is a tensor product of powers of several line bundles. That case is explored further in Section \ref{S4}. There, as a corollary of this theorem, we show that (\ref{Main-Convergence}) holds with fewer assumptions.

\begin{Theorem}\label{T1.3}
Let $(X,\omega)$ and $(L_p,h_p)$ be as in $(A)$ and $(B)$. If every $x\in X$ has a neighborhood $U$, with local frames $e_p$, such that the family $\{\phi_p/A_p\}$ is equicontinuous, and $\{(\phi_p-\rho_p)/A_p\}$ is uniformly bounded, then (\ref{Main-Convergence}) holds.
\end{Theorem}

\ii When $h_p$ satisfies (\ref{Defn ap}), then we may assume $h_p=g_p$, and so $\{(\phi_p-\rho_p)/A_p=0\}$ is automatically bounded. In general, the difference $(\phi_p-\rho_p)$ defines a function on all $X$ (see Section \ref{S2.2}). In the case where $(L_p,g_p)=(L^p,g^p)$ and $e_p=e^{\otimes p}$, we can take $A_p=p$, and $\{(\phi_p-\rho_p)/p\}$ is bounded whenever $\{\phi_p/p\}$ is. In this case, Theorem \ref{T1.5} shows that $\gamma_p/p$ converges.

\ii This leads us to the question of under what conditions does $\gamma_p/A_p$ converge? Due to Coman, Ma, and Marinescu, we know if $h_p=g_p$ and $c_1(L_p,h_p)/A_p$ converges, then so does $\gamma_p/A_p$. In the setting of our previous theorems we also know that if $c_1(L_p,h^{eq}_p)/A_p$ converges, then so does $\gamma_p/A_p$. So to give a partial answer to our question, we assume $c_1(L_p,h_p)/A_p$ converges, and ask under what conditions does $c_1(L_p,h^{eq})/A_p$ converge as well? We will state a theorem with sufficient conditions for such convergence. Before we do, we introduce a necessary proposition.

\begin{Proposition}\label{Prop}Let $X, V_\alpha$, and $\phi^\alpha$ be as in $(A)$ and $(B')$, then for all $\alpha,\beta$ there exist pluriharmonic functions $\psi^{\alpha\beta}$ on $V_\alpha\cap V_\beta$ such that
\[\phi^\alpha=\phi^\beta+\psi^{\alpha\beta}.\]
\end{Proposition}

\ii In this case, there exists a real closed $(1,1)$-current $T$ on $X$ defined by 

\begin{equation}\label{def of T current}
T|_{V_\alpha}=dd^c\phi^\alpha,
\end{equation}
where $d^c= \frac{1}{2\pi i}(\partial-\overline{\partial})$. Let $\{T\}$ denote the cohomology class of $T$, and fix a smooth form $\theta\in \{T\}.$ From the definition, we have 
\[
T=\theta+dd^c \varphi,
\]
where $\varphi\in L^1(X).$ Since $\phi^\alpha$ is continuous, it follows that $\varphi$ is continuous.

\begin{Theorem}\label{T1.4} Let $(X,\omega)$, $(L_p,h_p)$, and $A_p$ be as in $(A)$ and $(B')$. Suppose the following conditions hold:

\begin{enumerate}[label=(\alph*)]

 \item There exist $\delta_p\geq 0$ with $\delta_p\to0$ such that
\smallskip
\[T-\delta_p \omega \leq \frac{c_1(L_p,h_p)}{A_p} \text{ on }V_\alpha\text{ for all }\alpha.\]

\item There exist $\varrho\in PSH(X,\theta)$ and some $c>0$ such that 
\[\theta+dd^c\varrho\geq c\omega.\]

\end{enumerate}

Then 
\begin{equation}\label{T1.4-convergence}
\frac{\gamma_p}{A_p}\to \theta+dd^c\varphi^{eq}\text{ and }\ \frac{c_1(L_p,h_p^{eq})}{A_p}\to \theta+dd^c\varphi^{eq}
\end{equation}

weakly as currents, where

\[\varphi^{eq}:=\sup\{\psi\in PSH(X,\theta)\divides \psi\leq\varphi \}.\]
Moreover, 
\begin{equation}\label{alternate definition of varphi eq}
\varphi^{eq}=\left[\limsup_p\left(\frac{\varphi_p^{eq}-\varphi_p}{A_p}\right)\right]^*+\varphi,
\end{equation}
where $\varphi_p$ and $\varphi_p^{eq}$ are the global weights of $h_p$ and $h_p^{eq}$ (see section \ref{S2.3}).
\end{Theorem}

\ii Here $PSH(X,\theta)$ denotes the class of $\theta-$plurisubharmonic functions (see section \ref{S2.1})  Note that (\ref{T1.4-convergence}) and condition (b) are both independent of our choice of $\theta.$ The independence of (\ref{T1.4-convergence}) is shown in the proof of Theorem \ref{T1.4}. In the case of (b), if $\widetilde{\theta}\in \{T\}$, then $\theta$ and $\widetilde{\theta}$ are in the same cohomology class, so clearly there exists $\widetilde{\varrho}\in L^1(X)$ with $\widetilde{\theta}+dd^c\widetilde{\varrho}=\theta+dd^c\varrho.$

\ii We said we were interested in cases where $c_1(L_p,h_p)/A_p$ converges, and this is one of them, as if $\phi_p^\alpha/A_p\to \phi^\alpha$ uniformly, then

\[\frac{c_1(L_p,h_p)}{A_p}\to T\text{ weakly as currents, where }T\text{ is as in (\ref{def of T current}}).\]

\ii The function $\varphi^{eq}_p$ is defined with respect to a smooth metric $h_p^{0}$ in section \ref{S2.2}. One can show using the $C^\infty$ version of the first Cousin Problem, that for particular choices of $h_p^0$ and $\theta$, equation (\ref{alternate definition of varphi eq}) reduces to  $\varphi^{eq}=[\limsup_p (\varphi^{eq}_p/A_p)]^*$.

\ii In Lemma \ref{3.2.2} we will show directly that there exists $\psi\in PSH(X,\theta)$ with $ \psi\leq \varphi$, hence $\varphi^{eq}\not=\infty$. 

\ii The purpose of conditions (a) and (b) is to allow us to apply Demailly's $L^2$ estimates for $\overline{\partial}$, which requires some degree of positivity. 

\ii In Section \ref{S4} we will look at the case where there exists a continuous $(1,1)-$form $\Phi$ such that ${c_1(L_p,h_p)/A_p}\to \Phi$ uniformly. This assumption implies the convergence in $(B')$. It is examined in Corollaries \ref{C1.4.1} and \ref{C1.4.2}. If we make  the additional assumption $L_p=L^p$, then (b) is automatically satisfied whenever $L$ is big, a fact which can be used to prove a special case of Theorem \ref{T1.5}.

\ii This paper has the following structure. In Section \ref{S2}, we discuss the necessary information about quasi-plurisubharmonic functions, define global weights for our metrics, and give a construction of $h_p^{eq}$. Section \ref{S3} is devoted to proofs of the theorems above. In Section \ref{S4}, we present applications of our main results.

\section{Preliminaries}\label{S2}

\subsection{Quasi-plurisubharmonic functions}\label{S2.1} Recall that a quasi-plurisubharmonic (qpsh) function on $X$ is a function, $f:X\to [-\infty,\infty),$ that is locally the sum of a plurisubharmonic (psh) function and a smooth function. Given a real closed smooth $(1,1)-$form $\theta$, a $\theta-$plurisubharmonic ($\theta-$psh) function is a qpsh function $f:X\to[-\infty,\infty)$ such that $\theta+dd^cf\geq 0$.  We denote by $PSH(X,\theta)$ the class of all $\theta-$psh functions.

\ii In the setting of Theorem \ref{T1.1}, note the condition  $-M\omega \leq c_1(L_p,h_p)/A_p\leq M\omega$ implies $-M\omega\leq dd^c\phi_p/A_p\leq M\omega$ on any contractible Stein open set $V$ where $dd^c\phi_p= c_1(L_p,h_p).$ In this case $\pm\phi_p/A_p\in PSH (V, M\omega).$ This is one of the reasons we are interested in the notion of qpsh functions. Another is that they are necessary for us to define the equilibrium metrics $h_p^{eq}$. 

\ii One property will we need to note about $\theta$-psh functions, is that two are equal almost everywhere, then they are equal. We will use this fact throughout this paper.

\subsection{Global weights and the equilibrium metric}\label{S2.2} We first observe that if $\Omega$ is a contractible Stein open set in $\mathbb{C}^n$, then a result of Oka tells us $H^1(\Omega,\mathcal{O}^*)\cong H^2(\Omega,\Z)$, hence line bundles over $\Omega$ are trivial (see \cite[p.201]{Huy05} and \cite{Oka39}). In particular, if $V\subset X$ is a coordinate polydisc, then $L_p$ is trivialized on $V$ for all $p\geq 1$.

\ii We fix a finite open cover $\{U_\alpha\}$ of $X$ by coordinate polydiscs, such that for all $\alpha$, there exists a coordinate polydisc $V_\alpha$ with $U_\alpha\Subset V_\alpha$. Fix local frames $e_p^\alpha$. Let $\phi_p^\alpha$ and $\rho_p^\alpha$ denote the local weights of $h_p$ and $g_p$ on $V_\alpha$. As $X$ is compact, we assume that the $V_\alpha$ are the neighborhoods referred to in our theorems.

\ii In order to define global weights for our given metrics, as well as give a definition of $h_p^{eq}$, we fix a smooth metric $h_p^0$ on $L_p$. Define $\xi_p^\alpha\in C^\infty(V_\alpha)$ to be the local weights of $h_p^0,$ which are given by
\smallskip
\begin{equation}\label{def of xi p}
h_p^0(e_p^\alpha,e_p^\alpha)=e^{-2\xi_p^\alpha}.
\end{equation}

The global weights of $h_p$ and $g_p$ are the functions $\varphi_p:X\to\R$ and $\varrho_p:X\to[-\infty,\infty)$, defined by
\smallskip
\begin{equation}\label{Local Weights}
h_p=h_p^0e^{-\varphi_p}\text{ and }g_p=h_p^0e^{-2\varrho_p}.
\end{equation}

Note that $(\phi_p^\alpha-\rho_p^\alpha)=(\varphi_p-\varrho_p)$, which is a global function. 

\ii In order to define $h_p^{eq}$, we first set

\[\theta_p=c_1(L_p,h_p^0).\]

Let $\varphi_p^{eq}:X\to[-\infty,\infty)$ be the $\theta_p-$psh upper envelope
\smallskip
\begin{equation}\label{definition of varphi eq}
\varphi_p^{eq}=\sup\{\psi\in PSH(X,\theta_p)\divides \psi\leq \varphi_p\}.
\end{equation}

Observe that $\varphi_p$ is continuous and $\varrho_p$ is bounded above (it is psh), so $\varrho_p-C\leq \varphi_p$ for some $C\geq 0$. This tells us
\smallskip
\[\varrho_p-C\in \sup\{\psi\in PSH(X,\theta_p)\divides \psi\leq \varphi_p\},\]

so $\varphi_p^{eq}\in PSH(X,\theta_p)$ and $\varphi^{eq}_p\leq \varphi_p$. We define 

\[
h_p^{eq}=h_p^0 e^{-2\varphi_p^{eq}}.
\]

The global weight of $h_p^{eq}$ is $\varphi_p^{eq}.$ A local construction of the equilibrium metric is given by Berman in \cite{Ber09}. His definition is clearly equivalent to ours.

\subsection{Fubini-Study currents and Bergman kernels}\label{S2.3} Let $H_{(2)}^0(X,L_p)$ denote the Bergman space of square integrable sections of $L_p$ relative to $h_p$ and $\omega$, that is,
\[
H_{(2)}^0(X,L_p,h_p)=\left\{S\in H^0(X,L_p)\divides \|S\|^2_{h_p}<\infty\right\},
\]

where
\[\|S\|_{h_p}^2:=\int_X |S|_{h_p}^2\frac{\omega^n}{n!}.\]

When the metric is clear, the notation is shortened to $H_{(2)}^0(X,L_p)$. This space will be endowed with the inner product 

\[\left<S_1,S_2\right>=\int_X h_p(S_1,S_2) \frac{\omega^n}{n!}.\]

\ii Let $P_p$ be the Bergman kernel function of the space $H_{(2)}^0(X,L_p).$ For all $p\geq 1$, a global definition is given by fixing an orthonormal basis, $\{S^p_j\}$, of $H_{(2)}^0(X,L_p)$. We then define $d(p):= \dim H_{(2)}^0(X,L_p)$ and

\[
P_p(x):= \sum_{j=1}^{d(p)} |S^p_j(x)|_{h_p}^2.
\]

 The following is a well-known variational characterization of the Bergman kernel, which will be useful in our work,

\[
P_p= \sup_{S\in H_{(2)}^0(X,L_p)} \left(\frac{|S|_{h_p}^2}{\|S\|^2_{h_p}}\right).
\]

\ii We recall that $\gamma_p$ was the Fubini-Study current of $H_{(2)}^0(X,L_p).$ To define it explicitly, we let $U\subset X$ be a contractible Stein open set. Let $s^p_j\in\mathcal{O}(U)$ be defined by $S^p_j=s^p_j e_p$. Then
\smallskip
\[
\gamma_p|_{U}:=\frac{1}{2}\ {dd^c \log\left(\sum_{j=1}^{d(p)}|s^p_j|^2\right)}.
\]

Note, that the equation above defines a global current on $X$. One way to see that is the following equivalent definition. We consider the Kodaira map $\Psi_p:X\dashedrightarrow \mathbb{CP}^{k-1}$ defined by
\smallskip
\[x\mapsto [s^p_1(x):\ldots: s^p_{d(p)}(x)],\ \forall x\in U\bs V(s_1^p,\ldots s_{d(p)}^p),\]

where $V(s_1^p,\ldots, s_k^p)$ is the analytic variety. It is well known that this map is independent of our choice of $U.$ We then define
\smallskip
\[\gamma_p= \Psi_p^*\omega_{FS},\]

where $\omega_{FS}$ is the Fubini-Study form on $\mathbb{CP}^{k-1}$. 

\ii Note, the following is a well-known identity which will be helpful in our work: 

\begin{equation}\label{Fubin-Study-identity}
\gamma_p=\frac{dd^c \log P_p}{2}+c_1(L_p,h_p).
\end{equation}

\section{Main Results}\label{S3}

\subsection{Proofs of Theorems \ref{T1.1}, \ref{T1.2}, and \ref{T1.3}} We use the notation and definitions introduced in Sections \ref{S1} and \ref{S2}. We start by stating and proving two lemmas necessary for our proof of Theorem \ref{T1.1}.

\begin{Lemma}\label{S2-Continuity-Lemma}
  If $V\subseteq \mathbb{C}^n$ is a polydisc and $v\in L^1(V)$ such that for some $M\geq 0$ we have
  \[-M\omega\leq dd^c v\leq M\omega,\]
  then there exists $\widetilde{v}\in C(V)$ with $\widetilde{v}=v$ a.e. and $dd^c\widetilde{v}=dd^cv.$
\end{Lemma}
\emph{Proof.} Let $\zeta$ be a smooth potential for $M\omega$ on $V$. Since $-M\omega\leq dd^c v\leq M\omega$ it follows that $\zeta+v$ and $\zeta-v$ are both equal a.e. to psh functions on $V.$ Let $u,w$ be psh functions with $\zeta+v=u$ and $\zeta-v=w$ a.e. Note that $2\zeta=u+w$ a.e. Since both sides are psh we have $2\zeta=u+w$ everywhere. Then, $2\zeta-u=w$, so as $2\zeta-u$ is lower semicontinuous, and $w$ is as well. We see that $w$ is both upper and lower semicontinuous, so it is continuous, and $v=w$ a.e. 
\eproof

\begin{Lemma}\label{S2-Eventually-Equicontinuous?}
Let $U,V\subset X$ be open coordinate polydiscs with $U\Subset V$. Assume $v_p\in C(V)$ and $v\in L^1(V)$ such that $v_p\to v$ in $L^1$,
\smallskip
\[
-M\omega\leq dd^c v_p\leq M\omega.
\]

Then, for all $\epsilon>0$ there exists $p_0=p_0(\epsilon)$ and $\delta=\delta(\epsilon)$ such that if $z\in U,$ $r<\delta,$ and $p>p_0$, then 
\smallskip
\[\sup_{B(z,r)}v_p(z)-\inf_{B(z,r)}v_p(z)<\epsilon.\]
\end{Lemma}

\emph{Proof.} By Lemma \ref{S2-Continuity-Lemma} we can and do assume $v$ is $C$. Let $\zeta$ be a smooth real potential for $M\omega$ on $V$. By \cite[Theorem 3.2.12]{Hor94}, if $\zeta$ is a real smooth potential for $M\omega$, then we can and do assume $\zeta\pm v\in PSH(V).$ 

\ii From Hartog's lemma (see \cite[Theorem 3.2.13]{Hor94}) applied to $v_p$ and $v$ it follows that for all $K\Subset V$ we have
\smallskip
\begin{equation}\label{S2-Hartogs}
\limsup_{p\to\infty}\left[\sup_K v_p\right]\leq \sup_K v.
\end{equation}

\ii Let $B_1,\ldots,B_m$ be a finite cover of $U$ by balls such that $2B_j\Subset V$ for all $j$, where $2B_j$ is the dilation of $B_j$ by a factor of $2.$ Using (\ref{S2-Hartogs}), since there are finitely many $B_j$, we can choose $p_0$ large enough so that 
\smallskip
\begin{equation}\label{S2-sup}
\sup_{2B_j}v_p\leq \sup_{2B_j} v+\epsilon
\end{equation}

for all $j$ and $p> p_0.$ Similarly, we may assume for such $p$ we have 
\smallskip
\[\sup_{2B_j}(-v_p)\leq \sup_{2B_j}(-v)+\epsilon,\]
thus
\smallskip
\begin{equation}\label{S2-inf}
-\inf_{2 B_j}v_p\leq -\inf_{2B_j}v+\epsilon.
\end{equation}
\smallskip
\ii Let $\delta>0$ be the minimum of the radii of the $B_j.$ Additionally, assume $r<\delta$ and $z\in U.$ Since the sets $B_1,\ldots,B_m$ cover $U$, we deduce $z\in B_j$ for some $j,$ therefore $B(z,r)\subset 2B_j,$ and if $p\geq p_0$, then by (\ref{S2-sup}) and (\ref{S2-inf}) we have

\[
\sup_{B(z,r)} v_p-\inf_{B(z,r)}v_p\leq \sup_{2B_j} v_p-\inf_{2B_j} v_p\leq \sup_{2B_j} v-\inf_{2B_j}v+2\epsilon.
\]

As $v$ is continuous, by letting the radii of $B_j$ go to zero we have our desired upper bound. 
\eproof

\emph{Proof of Theorem \ref{T1.1}.} Recall that $\{U_\alpha\}$ is a finite cover of $X$ by coordinate polydiscs. Moreover, for all $\alpha$ there exists a coordinate polydisc $V_\alpha$ such that $U_\alpha\subseteq V_\alpha$. 

\ii Notice that for $p\geq 1$, the sets $V_\alpha$ and local frames $e_\alpha$ satisfy the assumptions of our theorem. That is, the local weights $\phi_p^\alpha$ and $\rho_p^\alpha$ defined by 
\[h_p(e_p^\alpha,e_p^\alpha)=e^{-2\phi_p^\alpha}\text{ and }g_p(e_p^\alpha,e_p^\alpha)=e^{-2\rho_p^\alpha},\]

form a family $\{\phi_p^\alpha/A_p,\rho_p^\alpha/A_p\}$ which is uniformly bounded in $L^1(V_\alpha).$ Additionally, for each $p\geq 1$, we recall that $h_p^0$ is a smooth metric on $L_p$ with global weights $\xi_p^\alpha$ with respect to $e_p^\alpha$. The global weights $\varphi_p$ and $\varrho_p$ are then defined by
\[h_p=h_p^0 e^{-2\varphi_p}\text{ and }g_p=h_pe^{-2\varrho_p}.\]

Furthermore, the equilibrium metric $h_p^{eq}$ is given by $h_p^{eq}=h_p^0 e^{-2\varphi_p^{eq}}$.

\ii For a fixed $p$, let $\{S_j\}$ be an orthonormal basis of $H_{(2)}^0(X,L_p)$ represented locally by $S_j=s_j^\alpha e_p^\alpha$. Observe that by (\ref{Fubin-Study-identity}) we have

\begin{equation}\label{gamma-p curvature equation}
\frac{\gamma_p-c_1(L_p,h_p^{eq})}{A_p}=dd^c\left[\left(\frac{(\log P_p)/2+\varphi_p}{A_p}\right)-\frac{\varphi_p^{eq}}{A_p}\right].
\end{equation}

Define
\smallskip
\[\beta_p=\frac{(\log P_p)/2+\varphi_p}{A_p}.\]

It follows that
\smallskip
\[\frac{\gamma_p-c_1(L_p,h_p^{eq})}{A_p}=dd^c\left(\beta_p-\frac{\varphi_p^{eq}}{A_p}\right),\]

so if we show $\beta_p-\varphi_p^{eq}/A_p\to 0$ in $L^1$ we will prove our claim. This happens if every subsequence has a subsequence where the convergence holds, so it suffices to show a single subsequence where convergence holds. 

\ii Let $\zeta^\alpha\in C^\infty(V_\alpha,\mathbb{R})$ be a real smooth potential for $M\omega$. By (\ref{T1.1-Inequality}), it follows that $\phi_p^\alpha/A_p+\zeta^\alpha\in C(V_\alpha)$ is psh. For all $\alpha$ define 
\smallskip
\[\delta=\delta(\alpha)=\text{dist}(\partial V_\alpha, U_\alpha),\]

where dist is the distance induced by $\omega$ on $X$.

For any $\epsilon\geq 0$, let
\smallskip
\[U_\alpha^{\epsilon}=\{x\in V_\alpha\divides \text{dist}(x,U_\alpha)< \epsilon\}.\]

\ii By the subaveraging property of psh functions and the $L^1(V_\alpha)$ uniform boundedness of $\{\phi_p^\alpha/A_p\}$, it follows that $\{\phi_p^\alpha/A_p+\zeta^\alpha\}$ is uniformly upper bounded on ${U_\alpha^{\delta/2}}.$ The same logic tells us $\{\zeta^\alpha-\phi_p^\alpha/A_p\}$ is uniformly upper bounded on ${U_\alpha^{\delta/2}}$. Since $\zeta^\alpha$ is continuous, putting these two bounds together allows us to conclude that  $\{\phi_p^\alpha/A_p+\zeta^\alpha\}$ is a uniformly bounded family of psh functions on ${U_\alpha^{\delta/2}}.$ From \cite[Theorem 3.2.12]{Hor94}, there exists a psh function $\widehat{\phi}^\alpha$ and a subsequence $\{{\phi_{p_j}^\alpha}/{A_{p_j}}+\zeta^\alpha\}$ such that
\smallskip
\[{\phi_{p_j}^\alpha}/{A_{p_j}}+\zeta^\alpha\to \widehat{\phi}^\alpha \text{ in }L^1(U_\alpha^{\delta/2}).\]

Since there are finitely many $\alpha$, by passing to a further subsequence, if necessary, we may assume
\[{\phi_{p_j}^\alpha}/{A_{p_j}}+\zeta^\alpha\to \widehat{\phi}^\alpha \text{ in }L^1(U_\alpha^{\delta/2})\text{ and }a.e.\]

for all $\alpha$ simultaneously. This subsequence is the one where we will show convergence holds.

\ii Let $\epsilon>0$. By Lemma \ref{S2-Eventually-Equicontinuous?} it follows that we can choose $r=r(\alpha,\epsilon)\in\R$ satisfying 
\smallskip
\[\delta/8>r>0\]

small enough and $j_0>0$ large enough so that
\smallskip
\begin{equation}\label{T1.1-upper-bound-phi-p}
\sup_{B(z,r)}\left({\phi_{p_j}^\alpha}/{A_{p_j}}\right)-\inf_{B(z,r)}\left({\phi_{p_j}^\alpha}/{A_{p_j}}\right)<\epsilon
\end{equation}

for all $j>j_0$ and $z\in U_\alpha$. As there are finitely many $\alpha$, WLOG we may assume the above holds for all $\alpha$. 

\ii Let $x\in U_\alpha.$ By our choice of $r$, it follows that $V:=B(x,2r)\subset U_\alpha^{\delta/4}$. Let $U=B(x,r)$ and $z\in U$. Suppose $S\in H_{(2)}^p(X,L_p)$ and define $s^\alpha\in\mathcal{O}(V_\alpha)$ by $S=s^\alpha e_p^\alpha$ on $V_\alpha$.  Since $s^\alpha$ is holomorphic, it follows from subaveraging that

\[|S(z)|^2_p=|s^\alpha(z)|^2e^{-2\phi_p^\alpha(z)}\leq \frac{e^{-2\phi_p^\alpha(z)}}{\lambda(B(z,r))}\int_{B(z,r)}|s^\alpha|^2d\lambda\]

where $\lambda$ is the Lebesgue measure. It follows that there exists $C_1>0$ such that

\[|S(z)|^2_p\leq\frac{ C_1e^{-2\phi_p^\alpha(z)}}{\lambda(B(z,r))}\int_{B(z,r)}|s^\alpha|^2\frac{\omega^n}{n!}\]

\[\leq \frac{C_1e^{-2\phi_p^\alpha(z)+\max_{\overline{B}(z,r)}2\phi_p^\alpha}}{\lambda(B(z,r))}\int_{B(z,r)}|s^\alpha|^2\left(e^{-2\phi_p^\alpha}\right)\frac{\omega^n}{n!}\]

\[= \frac{\|S\|_{h_p}^2 C_1e^{-2\phi_p^\alpha(z)+\max_{\overline{B}(z,r)}2\phi_p^\alpha}}{\lambda(B(z,r))}.\]

By the variational characterization of the Bergman kernel, we have

\[P_p(z)\leq \frac{C_1e^{-2\phi_p^\alpha(z)+2\max_{\overline{B}(z,r)}\phi_p^\alpha}}{\lambda(B(z,r))},\]

therefore
\smallskip
\[\frac{\log P_{p_j}(z)}{2A_{p_j}}\leq \frac{\log C_1-2n\log r-n\log \pi+\log n!}{2A_{p_j}}+\max_{\overline{B}(z,r)}\left(\frac{\phi_{p_j}^\alpha}{A_{p_j}}\right)-\frac{\phi_{p_j}^\alpha(z)}{A_{p_j}}.\]

By (\ref{T1.1-upper-bound-phi-p}) and by enlarging $j_0$, if necessary, it follows that if $j\geq j_0,$ then on $B(x,r)$ we have
\smallskip
\[
\frac{\log P_{p_j}}{2A_{p_j}}\leq \frac{-2n\log r}{2A_{p_j}}+\epsilon.
\]
By covering $X$ with finitely many such $B(x,r)$, we may assume WLOG that the above inequality holds on $X$.

\ii Note that
\smallskip
\begin{equation}\label{T1.1-Proof-BD-By-phi}
\frac{\log P_{p_j}}{2}+{\varphi_{p_j}}-A_{p_j}\epsilon +\frac{2n\log r}{2}\leq{\varphi_{p_j}},
\end{equation}

and
\smallskip
\[dd^c\left(\frac{\log P_{p_j}}{2}+{\varphi_{p_j}}-\epsilon A_{p_j}+\frac{2n\log r}{2}\right)=\gamma_{p_j}-\theta_{p_j},\]

so the L.H.S. of (\ref{T1.1-Proof-BD-By-phi}) is $\theta_{p_j}-$psh. By the definition of $\varphi_p^{eq}$, we have
\smallskip
\[\frac{\log P_{p_j}}{2}+{\varphi_{p_j}}-\epsilon A_{p_j}+\frac{2n\log r}{2}\leq\varphi_{p_j}^{eq},\]

so 
\smallskip
\[\frac{\log P_{p_j}}{2A_{p_j}}+\frac{\varphi_{p_j}}{A_{p_j}}\leq \frac{-2n\log r}{2A_{p_j}}+\epsilon+\frac{\varphi^{eq}_{p_j}}{A_{p_j}},\]

and
\smallskip
\begin{equation}\label{Pro24}
\beta_{p_j}-\frac{\varphi_{p_j}^{eq}}{A_{p_j}}\leq \frac{-2n\log r}{2A_{p_j}}+\epsilon.
\end{equation}

This is the upper bound we shall use.

\ii To show the lower bound, we first define $t_p=1/\sqrt{a_p}$ and set

\[\widetilde{h}_p=h_p^{0}e^{-2(1-t_p)\varphi_{p}^{eq}-2t_p\varrho_p}.\]

We compute 
\smallskip
\[c_1(L_p,\widetilde{h_p})=\theta_p+(1-t_p)dd^c\varphi_p^{eq}+t_pdd^c\varrho_p\]

\[=(1-t_p)\left(\theta_p+dd^c \varphi_p^{eq}\right)+t_p\left(\theta_p+dd^c\varrho_p\right)\geq t_p\left(\theta_p+dd^c\varrho_p\right)\]

\[\geq \frac{a_p}{\sqrt{a_p}} \omega=\sqrt{a_p}\omega.\]

Here, $\sqrt{a_p}\to\infty$ as $p\to\infty$. This shows that $\widetilde{h}_p$ satisfies the assumptions of Demailly's $L^2$ estimates for $\overline{\partial}$ (see \cite[Theorem 2.5]{CMM17}). As in the proof of \cite[Theorem 1.1]{CMM17}, for all $p$ large enough, we use the Ohsawa-Takegoshi Extension Theorem and Demailly's estimates for $\overline{\partial}$ (see \cite{OT87} and \cite{Dem82}, resp.) to show the following. There exists $C_2>0$ such that for all $z\in U_\alpha$ with $\rho^\alpha_p(z)\not=-\infty$, there exists $S_{z,p}\in H^0(X,L_p)$ such that $S_{z,p}(z)\not=0$ and
\smallskip
\begin{equation}\label{Pro29}
\|S_{z,p}\|_{\widetilde{h}_p}^2\leq C_2|S_{z,p}(z)|_{\widetilde{h}_p}^2.
\end{equation}

\ii We compute
\smallskip
\[\widetilde{h}_p=h_{p}^0e^{-2(1-t_p)\varphi_p^{eq}-2t_p\varrho_p}=h_pe^{2\varphi_p-2(1-t_p)\varphi_p^{eq}-2t_p\varrho_p}\]
\[=h_pe^{2(1-t_p)(\varphi_p-\varphi_p^{eq})}e^{2t_p(\varphi_p-\varrho_p)}\geq h_pe^{2t_p(\varphi_p-\varrho_p)}= h_pe^{2t_p\left(\phi_p^\alpha-\rho_p^\alpha\right)}.\]

Since $\{\rho_p^\alpha/A_p\}$ is a family of psh functions uniformly bounded in $L^1(V_\alpha)$ it is locally uniformly bounded above. Then, as $\{\phi_p^\alpha/A_p\}$ is locally uniformly bounded, it follows that $\{(\phi_p^\alpha-\rho_p^\alpha)/A_p\}$ is bounded below on $U_\alpha$ for all $\alpha.$ Since $(\phi_p^\alpha-\rho_p^\alpha)/A_p$ forms a global function, there exists $D\in\R$ with 
\smallskip
\[(\phi_p^\alpha-\rho_p^\alpha)/A_p\geq D,\]

for all $p,\alpha$. Therefore,
\smallskip
\begin{equation}\label{Pro30}
\widetilde{h}_p\geq h_pe^{2t_p\left(\phi_p^\alpha-\rho_p^\alpha\right)}\geq h_pe^{2t_pDA_p}.
\end{equation}

From above and (\ref{Pro29}), we deduce that $S_{z,p}\in H_{(2)}^0(X,L_p)$ for all $p$ large enough. 

\ii From (\ref{Pro30}) and (\ref{Pro29}) we obtain
\smallskip
\[\|S_{z,p}\|_{h_p}^2e^{2t_pDA_p}\leq C_2|S_{z,p}(z)|_{h_p}^2\exp\left[{2\varphi_p(z)-2(1-t_p)\varphi_p^{eq}(z)-2t_p\varrho_p(z)}\right],\]

so
\smallskip
\[\frac{\exp\left[2{t_pDA_p-2\varphi_p(z)+2(1-t_p)\varphi_p^{eq}(z)+2t_p\varrho_p(z)}\right]}{C_2}\leq \frac{|S_{z,p}(z)|_{h_p}^2}{\|S_{z,p}\|_{h_p}^2}\leq P_p(z).\]

We compute
\[t_pDA_p+(1-t_p)\varphi_p^{eq}(z)+t_p\varrho_p(z)-\frac{\log(C_2)}{2}\leq \frac{\log P_p(z)}{2}+\varphi_p(z)\]

\[t_pD+\frac{-t_p\varphi_p^{eq}(z)+t_p\varrho_p(z)-\log(C_2)/2}{A_p}\leq \left(\beta_p-\frac{\varphi_p^{eq}}{A_p}\right)(z).\]

As $\varphi_p\geq \varphi_p^{eq}$ we get
\smallskip
\[t_pD+t_p\frac{(\varrho_p-\varphi_p)(z)}{A_p}-\frac{\log C_2}{2A_p}\leq \left(\beta_p-\frac{\varphi_p^{eq}}{A_p}\right)(z).\]

This serves as our lower bound. 

\ii Combining the above inequality with that of (\ref{Pro24}), it follows that for $j$ large enough, we have

\[t_pD+t_p\frac{(\varrho_p-\varphi_p)(z)}{A_p}-\frac{\log C_2}{2A_p}\leq \left(\beta_{p_j}-\frac{\varphi_{p_j}^{eq}}{A_{p_j}}\right)(z)\leq -\frac{2n\log r}{2A_{p_j}}+\epsilon\]

for a.e. $z\in U$. By compactness it holds for a.e. $z\in X.$ Then, 

\[\left|\beta_{p_j}-\frac{\varphi_{p_j}^{eq}}{A_{p_j}}\right|(z)\leq \left|t_pD+t_p\frac{(\varrho_p-\varphi_p)(z)}{A_p}-\frac{\log C_2}{2A_p}\right|+\left|\frac{2n\log r}{2A_{p_j}}\right|+\epsilon,\]

Recall that $(\phi_p^\alpha-\rho_p^\alpha)=(\varphi_p-\varrho_p)|_{V_\alpha}$ is bounded in $L^1(U_\alpha)$ for all $\alpha$. So, $(\varphi_p-\varrho_p)$ is bounded in $L^1(X).$ As well, $t_p\to 0$, so by integrating and letting $j\to\infty$, we find

\[0\leq \limsup_j
\int_X\left|\beta_{p_j}-\frac{\varphi_{p_j}^{eq}}{A_{p_j}}\right|\frac{\omega^n}{n!}\leq\epsilon\int_U\frac{\omega^n}{n!}.\]

By letting $\epsilon\to 0$ we see that
\smallskip
\[\beta_{p_j}-\frac{\varphi_{p_j}^{eq}}{A_{p_j}}\to 0\]

in $L^1(X)$. 
\eproof

\emph{Proof of Theorem \ref{T1.2}} Like in Theorem \ref{T1.1}, the family $\{U_\alpha\}$ forms a finite cover of $X$ by coordinate polydiscs, $V_\alpha$ is a coordinate polydisc containing $U_\alpha$, and $e_p^\alpha$ is a local frame of $L_p$ on $V_\alpha$. Here the global weights $\phi_p^\alpha$ and $\rho_p^\alpha$ of $h_p$ and $g_p$ satisfy the property that $\{\phi_p^\alpha/A_p\}$ is uniformly bounded and equicontinuous, and $\{\rho_p^\alpha/A_p\}$ is uniformly bounded in $L^1$.

\ii Define
\smallskip
\[\beta_p=\frac{\left(\log P_p\right)/{2}+\varphi_p}{A_p}.\]

As we established in the preceding proof, if we show $\beta_p-\varphi_p^{eq}/A_p\to 0$ in $L^1$, then the proof is completed. 

\ii Let $\epsilon>0$ and $x\in X$. Using our subaveraging argument from the proof of Theorem \ref{T1.1} and the equicontinuity argument in the proof of \cite[Proposition 4.4]{CMM17}, it follows that 

\begin{equation}\label{9-6-24 Equation}A_p\beta_p-\frac{\log C_1-2n\log r}{2}-A_p\epsilon\leq  {\varphi_p},
\end{equation}

for some $C_1\in\R$. The L.H.S. is $\theta_p-$psh, hence by our arguments in the proof of Theorem $1.1$, the upper bound $\varphi_p$ may be replaced by $\varphi_p^{eq}$. Therefore, 
\[\left(\beta_p-\frac{\varphi_p^{eq}}{A_p}\right)\leq \frac{\log C_1-2n\log r}{2A_p}+\epsilon.\]
This serves as our upper bound. 

\ii For the lower bound, we apply the same argument as in the proof of Theorem \ref{T1.1}, which only required $\{\phi_p^\alpha/A_p\}$ to be locally uniformly bounded and $\{\rho_p^\alpha/A_p\}$ to be bounded in $L^1$. These methods allow us to deduce that for some neighborhood $U$ of $x$, if $z\in U$, then 
\smallskip
\[t_{p}C_2+\frac{t_{p}(\varrho_p-\varphi_p)(z)-\log(C_3)/2}{A_{p}}\leq \left(\beta_p-\frac{\varphi_p^{eq}}{A_p}\right)(z)\leq \frac{\log C_1-2n\log r}{2A_p}+\epsilon\]

for some $C_2,C_3>0$. This is the same situation we had at the end of our proof of Theorem \ref{T1.1}, so we deduce $\beta_p-\varphi_p^{eq}/A_p\to 0$ in $L^1$.
\eproof

\emph{Proof of Theorem \ref{T1.3}.} Here $U_\alpha,V_\alpha, e_p^\alpha$ and the global weights $\phi_p^\alpha,\rho_p^\alpha$ are defined as in the previous two theorems. The conditions imposed on the weights by our assumptions are that $\{\phi_p/A_p\}$ is equicontinuous and $\{(\phi_p-\rho_p)/A_p\}$ is uniformly bounded. We define $\beta_p$ as in the previous two proofs. As before, it suffices to show $\beta_p-\varphi_p^{eq}/A_p\to 0$ in $L^1.$

\ii As in the proof of Theorem \ref{T1.2}, (\ref{9-6-24 Equation}) holds. For the lower bound, we first fix $z\in X$, then proceed as in our proof of Theorem \ref{T1.1}. Set $t_p=1/\sqrt{A_p}$ and 
\smallskip
\[\widetilde{h}_p=h_p^{0}e^{-2(1-t_p)\varphi_{p}^{eq}-2t_p\varrho_p}.\]

This definition is the same as in Theorem \ref{T1.1}, hence it satisfies

\[c_1(L_p,\widetilde{h_p})\geq \sqrt{a_p}\omega.\]

In exactly the same fashion as in our proof of Theorem \ref{T1.1}, for $p$ large enough, we conclude there exists $C_2>0$ such that for all $z\in U_\alpha$ with $\rho^\alpha_p(z)\not=-\infty$, there exists $S_{z,p}\in H^0(X,L_p)$ such that $S_{z,p}(z)\not=0$ and
\smallskip
\begin{equation}\label{Pro14}
\|S_{z,p}\|_{\widetilde{h}_p}^2\leq C_1|S_{z,p}|_{\widetilde{h}_p}^2(z).
\end{equation}

\ii As we have computed in our previous proofs, we have
\smallskip
\[\widetilde{h}_p\geq h_pe^{2t_p\varphi_p-2t_p\varrho_p}.\]

Since $\{(\varphi_p-\varrho_p)/A_p\}$ is uniformly bounded, it follows that there exists $D$ independent of $p$ such that
\smallskip
\[\widetilde{h}_p\geq h_pe^{2A_pt_pD},\]
and so $S_{z,p}\in H^0_{(2)}(X,L_p,h_p)$. From this and $(\ref{Pro14})$, we deduce
\smallskip
\[\frac{\exp\left[{2A_pt_pD-2\varphi_p(z)+2(1-t_p)\varphi_p^{eq}(z)+2t_p\varrho_p(z)}\right]}{C_1}\leq P_p(z).\]

Here $z\in X$ was arbitrary, so the above holds on all $X$, and we deduce
\smallskip
    \[t_pD+\frac{-t_p\varphi_p^{eq}+t_p\varrho_p-\log C_1}{A_p}\leq \beta_p-\frac{\varphi_p^{eq}}{A_p}.\]

Since $\varphi_p\geq \varphi_p^{eq}$, we find
\smallskip
\[t_pD+\frac{-t_p\varphi_p+t_p\varrho_p-\log C_1}{A_p}\leq \beta_p-\frac{\varphi_p^{eq}}{A_p}.\]
By the uniform boundedness of $\{(\varphi_p-\varrho_p)/A_p\}$, it follows that there exists $D_2>0$ such that
\smallskip
\[t_pD_2+\frac{-\log C_1}{A_p}\leq \beta_p-\frac{\varphi_p^{eq}}{A_p}.\]

Using our lower bound and (\ref{9-6-24 Equation}), one concludes as in the proof of Theorem \ref{T1.1}
 that
 \[\beta_p-\frac{\varphi_p^{eq}}{A_p}\to 0\text{ in }L^1(X).\]

\eproof

\subsection{Proofs of Proposition \ref{Prop} and Theorems \ref{T1.4}} We start with the proof of Proposition \ref{Prop}.

\emph{Proof.} Let $k_p^{\alpha\beta}$ be the transition functions for $L_p$ with respect to the open cover $\{V_\alpha\}$. Set $\psi_p^{\alpha\beta}=\log|k^{\alpha\beta}_p|$.
Then, we have the formula $\psi_p^{\alpha\beta}=\phi_p^\alpha-\phi_p^\beta$. Since $\phi_p^\alpha/A_p$ converges uniformly for all $p$, so does $\psi_p^{\alpha\beta}/A_p.$ Let the limit be $\psi^{\alpha\beta}.$ It follows that 
\begin{equation}\label{Cousin-Like-EQ}
\psi^{\alpha\beta}=\lim_p \left(\frac{\phi_p^\alpha-\phi_p^\beta}{A_p}\right)=\phi^\alpha-\phi^\beta.
\end{equation}
Moreover, $k_{p}^{\alpha\beta}$ is holomorphic and non-vanishing, hence  $\psi_p^{\alpha\beta}=e^{\log|k_p^{\alpha\beta}|}$ is pluriharmonic, and so is the uniform limit $\psi^{\alpha\beta}=\lim_p (\psi_p^{\alpha\beta}/A_p)$. This shows that $\psi^{\alpha\beta}$ satisfies the desired conclusions.
\eproof

\ii The following lemmas will be used to prove Theorem \ref{T1.4}. 

\begin{Lemma}
If $(X,\omega)$ and $(L_p,h_p)$ satisfy $(A)$ and $(B')$, and $\psi^{\alpha\beta}$ are functions such that $\phi^\alpha=\phi^\beta+ \psi^{\alpha\beta},$
then there exists $\xi^\alpha\in C^\infty$ with 
\[\xi^\alpha=\xi^\beta+\psi^{\alpha\beta}.\]
\end{Lemma}

\emph{Proof.} Note that (\ref{Cousin-Like-EQ}) implies 
\[\psi^{\alpha\beta}=-\psi^{\beta\alpha}\text{ and }\psi^{\alpha\beta}+\psi^{\beta\gamma}+\psi^{\gamma\alpha}\text{ for all }\alpha,\beta,\gamma.\]
By the $C^\infty$ version of the first Cousin Problem \cite[Proposition 6.1.7]{Kra92} there exist functions $\xi^\alpha\in C^\infty(V_\alpha)$ for all $\alpha$ such that
\[\xi^\alpha=\xi^\beta+\psi^{\alpha\beta}\text{ for all }\alpha,\beta.\]
This is what we desired to show.

\eproof

\ii Recall that the local weights $\varrho_p,$ $\varphi_p$ of $g_p$ and $h_p$ were defined in Section \ref{S2.2} (see (\ref{Local Weights})). As were $\varphi_p^{eq}$ and $\xi_p^\alpha$ (see (\ref{definition of varphi eq}) and (\ref{def of xi p})). We also recall that $\varphi_p=\xi_p^\alpha+\phi_p^\alpha$. In general, similar relationships exist for the local and global weights of any metric.

\ii In the following lemmas and the proof of Theorem \ref{T1.4}, we let $\theta$ be the $(1,1)-$form such that $\theta|_{V_\alpha}=dd^c\xi^\alpha$, where $\xi^\alpha$ are as in the previous lemma. Also, we will set $\varphi =\phi^\alpha-\xi^\alpha$. 

\begin{Lemma}\label{3.2.2} Let $(X,\omega)$ and $(L_p,h_p)$ satisfy $(A)$ and $(B')$, and suppose the following conditions hold.

\begin{enumerate}[label=(\alph*)]

 \item There exist $\delta_p\geq 0$ with $\delta_p\to0$ such that
\smallskip
\[dd^c \phi^\alpha -\delta_p \omega \leq \frac{c_1(L_p,h_p)}{A_p}\text{ on }V_\alpha.\]

\item There exist $\varrho\in PSH(X,\theta)$ and some $c>0$ such that 
\[\theta+dd^c\varrho\geq c\omega.\]

\end{enumerate}

Then there exists a subsequence $\{(\varphi_{p_j}^{eq}-\varphi_{p_j})/A_{p_j}\}$ such that
\[\left(\frac{\varphi_{p_j}^{eq}-\varphi_{p_j}}{A_{p_j}}\right)\to\widetilde{\varphi^{eq}}-\varphi\text{ in }L^1(X),\]
where
\begin{equation}\label{widetilde varphi eq is in the desired set}
\widetilde{\varphi^{eq}}:=\left[\limsup_{j\to\infty}\left(\frac{\varphi_{p_j}^{eq}-\varphi_{p_j}}{A_{p_j}}\right)\right]^*+\varphi\in\{\psi\in PSH(X,\theta)\divides \psi\leq \varphi\}.
\end{equation}
\end{Lemma}

\emph{Proof.} First, since $\varphi_p$ is continuous we may assume, by subtracting a constant, that $A_p(\varrho-\varphi)\leq 0$. We compute
\[dd^c(A_p(\varrho-\varphi)+\phi_p^\alpha)=A_p\left(dd^c\varrho+\theta\right)+dd^c(\phi_p^\alpha-A_p\phi^\alpha)\geq A_p c\omega -A_p \delta_p \omega,\]

hence for $p$ large enough it follows that $[A_p(\varrho-\varphi)+\phi_p^\alpha]$ is psh. 

\ii Let $p$ be large enough such that $[A_p(\varrho-\varphi)+\phi_p^\alpha]$ is psh. Since $A_p(\varrho-\varphi)+\varphi_p\in PSH(X,\theta_p)$ and $A_p(\varrho-\varphi)+\varphi_p\leq \varphi_p,$ it follows from the definition of $\varphi_p^{eq}$ that $A_p(\varrho-\varphi)+\phi^\alpha_p\leq\xi_p^{\alpha}+\varphi_p^{eq}\leq \phi_p^\alpha,$ therefore 
\[(\varrho-\varphi)+\phi^\alpha_p/A_p\leq(\xi_p^{\alpha}+\varphi_p^{eq})/A_p\leq \phi_p^\alpha/A_p.\]
By the local uniform convergence $\phi_p^\alpha/A_p\to\phi^\alpha$, the L.H.S. and R.H.S. of the previous equation are uniformly bounded in $L^1_{loc}(V_\alpha)$, thus $\{(\xi_{p}^{\alpha}+\varphi_{p}^{eq})/A_{p}\}$ is as well. As a consequence, there exists a subsequence $\{(\xi_{p_j}^{\alpha}+\varphi_{p_j}^{eq})/A_{p_j}\}$ which converges in $L^1(U_\alpha)$ and a.e. to a function $\widetilde{\phi}^{eq,\alpha}\in PSH (V_\alpha)$ (see \cite[Theorem 3.2.12]{Hor94}). As there are finitely many $\alpha$, we may assume the convergence holds for all $\alpha$ simultaneously. 

\ii We claim that
\begin{equation}\label{almost-tilde-phi-eq}
  \widetilde{\phi}^{eq,\alpha}-\phi^{\alpha}=\left[\limsup_{j\to\infty} \left(\frac{\varphi_{p_j}^{eq}-\varphi_{p_j}}{A_{p_j}}\right)\right]^*,
\end{equation}

To prove the claim, note that from \cite[Theorem 2.6.3]{Kli91} we can deduce

\[\widetilde{\phi}^{eq,\alpha}=\left[\limsup_{j\to\infty}\left(\frac{\xi_{p_j}^\alpha+\varphi^{eq}_{p_j}}{A_{p_j}}\right)\right]^*.\]

Due to the uniform convergence $\phi^\alpha_p/A_p\to \phi^\alpha$ and the continuity of $\phi^\alpha$, it follows that
\smallskip
\[\widetilde{\phi}^{eq,\alpha}-\phi^{\alpha}=\left[\limsup_{j\to\infty}\left(\frac{\xi_{p_j}^\alpha+\varphi^{eq}_{p_j}}{A_{p_j}}\right)\right]^*-\phi^\alpha= \left[\limsup_{j\to\infty}\left(\frac{\xi_{p_j}^\alpha+\varphi^{eq}_{p_j}}{A_{p_j}}-\phi^\alpha\right)\right]^*\]

\[=\left[\limsup_{j\to\infty}\left(\frac{\xi_{p_j}^\alpha+\varphi^{eq}_{p_j}-\phi_{p_j}^\alpha}{A_{p_j}}\right)\right]^*= \left[\limsup_{j\to\infty} \left(\frac{\varphi_{p_j}^{eq}-\varphi_{p_j}}{A_{p_j}}\right)\right]^*.\]
Which gives us (\ref{almost-tilde-phi-eq}).

\ii By definition
\[\widetilde{\varphi^{eq}}=\widetilde{\phi}^{eq,\alpha}-\xi^\alpha\]
is $\theta-$psh. Moreover, by (\ref{almost-tilde-phi-eq}), since $\varphi_p^{eq}\leq \varphi_p$ for all $p$, we have $\widetilde{\varphi^{eq}}\leq \varphi$ which gives us (\ref{widetilde varphi eq is in the desired set}). To finish off the computation we note that
\[\frac{\varphi_{p_j}^{eq}-\varphi_{p_j}}{A_{p_j}}=\frac{(\xi_{p_j}^\alpha+\varphi_{p_j}^{eq})-\phi^{\alpha}_{p_j}}{A_{p_j}}\to \widetilde{\phi}^{eq,\alpha}-\phi^{\alpha}=\widetilde{\varphi^{eq}}-\varphi,\]
where the convergence is in $L^1(U_\alpha)$. Since this holds on all $U_\alpha,$ it follows that $(\varphi_{p_j}^{eq}-\varphi_{p_j})/A_{p_j}\to \widetilde{\varphi^{eq}}-\varphi$ in $L^1(X)$. This completes the proof.
\eproof

\ii For the following proof we will use the definition of $\widetilde{\varphi^{eq}}$ from the previous Lemma.\\

\emph{Proof of Theorem \ref{T1.4}.} Let $\epsilon>0$ and $x\in X$. Although we are working with currents, it still suffices to show the claim holds for a subsequence. 

Since $\{\phi_p^\alpha/A_p\}$ is uniformly convergent on $V_\alpha$ it is equicontinuous on $U_\alpha$, hence by our arguments in the proof of Theorem \ref{T1.2}, (\ref{9-6-24 Equation}) holds. That is, there exists a $C_1\in\R$ such that
\begin{equation}\label{you need at end of this proof}
\frac{\log P_{p}}{2A_{p}}-\left(\frac{\varphi_{p}^{eq}-\varphi_{p}}{A_{p}}\right)\leq \frac{\log C_1-2n\log r}{2A_{p}}+\epsilon.
\end{equation}

From the above and Lemma \ref{3.2.2}, it follows that there exists a subsequence $\{p_j\}$ such that for all $j$, we have

\[
\frac{\log P_{p_j}}{2A_{p_j}}-(\widetilde{\varphi^{eq}}-\varphi)\leq \frac{\log C_1-2n\log r}{2A_{p_j}}+\epsilon+\left[\left(\frac{\varphi_{p_j}^{eq}-\varphi_{p_j}}{A_{p_j}}\right)-(\widetilde{\varphi^{eq}}-\varphi)\right].
\]

By the definition of $\varphi^{eq}$, it follows that

\begin{equation}\label{T1.4-EQ-3}
\frac{\log P_{p_j}}{2A_{p_j}}-(\varphi^{eq}-\varphi)\leq \frac{\log C_1-2n\log r}{2 A_{p_j}}+\epsilon+\left[\left(\frac{\varphi_{p_j}^{eq}-\varphi_{p_j}}{A_{p_j}}\right)-(\widetilde{\varphi^{eq}}-\varphi)\right].
\end{equation}

This shall serve as our upper bound. 

\ii By continuity of $\varphi_p$, we assume WLOG that $\varrho\in PSH(X,\theta)$ and $(\varrho-\varphi)\leq 0.$ From (b), we have
\[-\delta_p \omega \leq c_1(L_p,h_p)/A_p-dd^c \phi^\alpha=dd^c(\phi_p^\alpha/A_p-\phi^\alpha).\]

Define $t_{p_j}=\sqrt{\delta_{p_j}}+1/\sqrt{A_{p_j}},$ and set $\widetilde{h_{p_j}}= h_{p_j}e^{-2(1-t_{p_j})A_{p_j}(\varphi^{eq}-\varphi)-2t_{p_j}A_{p_j}(\varrho-\varphi)}.$ We compute

\[c_1(L_{p_j},\widetilde{h_{p_j}})=dd^c\phi_{p_j}^\alpha+(1-t_{p_j})A_{p_j}dd^c(\varphi^{eq}-\varphi)+t_{p_j} A_{p_j}dd^c(\varrho-\varphi)\]
\vspace{-1mm}
\[=dd^c(\phi_{p_j}^\alpha-A_{p_j}\phi^\alpha)+(1-t_{p_j})A_{p_j}dd^c(\varphi^{eq}+\xi^\alpha)+t_{p_j} A_{p_j}dd^c(\varrho+\xi^\alpha)\]
\vspace{-1.5mm}
\[\geq -\delta_{p_j}A_{p_j}\omega+\sqrt{\delta_{p_j}}A_{p_j}c\omega+c\sqrt{A_{p_j}}\omega\geq c\sqrt{A_{p_j}}\omega.\]

Since $c\sqrt{A_{p_j}}\to\infty$ as $j\to \infty$, $\widetilde{h_{p_j}}$ satisfies the assumptions of Demailly's $L^2$ estimates for $\overline{\partial}$. By our arguments from Theorem \ref{T1.1} using Ohsawa-Takegoshi extension and Demailly's estimate for $\overline{\partial}$, there is a $C_2>0$ such that for all $z\in U_\alpha$ with $\rho_p^\alpha(z)\not=-\infty$, there exists $S_{z,j}\in H^0(X,L_{p_j})$ satisfying $S_{z,j}(z)\not=0$ and
\smallskip
\begin{equation}\label{Sxj-inequality}
\|S_{z,j}\|^2_{\widetilde{h_{p_j}}}\leq C_2|S_{z,j}(x)|^2_{\widetilde{h_{p_j}}}.
\end{equation}

\ii Since $(\varphi^{eq}-\varphi)$ and $(\varrho-\varphi)$ are negative, we have $\widetilde{h_{p_j}}\geq h_{p_j},$ and so $S_{z,j}\in H^0_{(2)}(X,L_{p_j},h_{p_j})$. From this and (\ref{Sxj-inequality}) it follows that \[\|S_{z,j}\|_{{h_{p_j}}}\leq C_2|S_{z,j}(x)|_{{h_{p_j}}}e^{-2(1-t_{p_j})A_{p_j}(\varphi^{eq}-\varphi)(x)-2t_{p_j}A_{p_j}(\varrho-\varphi)(x)},\] thus
\[\frac{e^{2(1-t_{p_j})A_{p_j}(\varphi^{eq}-\varphi)(x)+2t_{p_j}A_{p_j}(\varrho-\varphi)(x)}}{C_2}\leq\frac{|S_{z,j}(x)|^2_{{h_{p_j}}}}{\|S_{z,j}\|^2_{{h_{p_j}}}}.\]

By the variational characterization of the Bergman kernel, it follows that
\[(1-t_{p_j})(\varphi^{eq}-\varphi)(x)+t_{p_j}(\varrho-\varphi)(x)-\frac{\log C_2}{2A_{p_j}}\leq \frac{\log P_{p_j}(x)}{2A_{p_j}},\]

hence
\[-t_{p_j}(\varphi^{eq}-\varphi)(x)+t_{p_j}(\varrho-\varphi)(x)-\frac{\log C_2}{2A_{p_j}}\leq \frac{\log P_{p_j}(x)}{2A_{p_j}}-(\varphi^{eq}-\varphi)(x).\]

As $(\varphi^{eq}-\varphi)$ is negative, we have
\[t_{p_j}(\varrho-\varphi)(x)-\frac{\log C_2}{2A_{p_j}}\leq \frac{\log P_{p_j}(x)}{2A_{p_j}}-(\varphi^{eq}-\varphi)(x).
\]

Combining with (\ref{T1.4-EQ-3}) yields
\begin{equation}\label{T1.5-another-equation}
t_{p_j}(\varrho-\varphi)(x)-\frac{\log C_2}{2A_{p_j}}\leq \frac{\log P_{p_j}(x)}{2A_{p_j}}-(\varphi^{eq}-\varphi)(x)
\end{equation}
\[\leq  \frac{\log C_1-2n\log r}{2 A_{p_j}}+\epsilon+\left[\left(\frac{\varphi_{p_j}^{eq}-\varphi_{p_j}}{A_{p_j}}\right)-(\widetilde{\varphi^{eq}}-\varphi)\right].\]

Note, by Lemma \ref{3.2.2} we have $(\varphi_{p_j}^{eq}-\varphi_{p_j})/A_{p_j}\to (\widetilde{\varphi^{eq}}-\varphi)$ in $L^1(X)$ and by definition $t_{p_j}\to 0$, so from computations similar to those in the proof of Theorem \ref{T1.1} we have

\smallskip
\[\limsup_j \left\|\frac{\log P_{p_j}}{2A_{p_j}}-(\varphi^{eq}-\varphi)\right\|_{L^1(X)}\leq \epsilon\int_X\frac{\omega^n}{n!}.\]

Letting $\epsilon\to 0$, we find that
\[\frac{\log P_{p_j}}{2A_{p_j}}-(\varphi^{eq}-\varphi)\to 0\]

in $L^1(X)$. Therefore 
\begin{equation}\label{T1.5-first-convergence}\frac{\gamma_{p_j}}{A_{p_j}}\to T+ dd^c\varphi^{eq}-dd^c\varphi=\theta+dd^c\varphi^{eq}
\end{equation}

weakly as currents. 

\ii Now, we address the second half of (\ref{T1.4-convergence}). By the definition of $\varphi^{eq}$ and Lemma \ref{3.2.2} we have 
\[\left(\frac{\varphi_{p_j}^{eq}-\varphi_{p_j}}{A_{p_j}}\right)-\left[\left(\frac{\varphi_{p_j}^{eq}-\varphi_{p_j}}{A_{p_j}}\right)-(\widetilde{\varphi^{eq}}-\varphi)\right]= (\widetilde{\varphi^{eq}}-\varphi)\leq  (\varphi^{eq}-\varphi).\]
With this and (\ref{T1.5-another-equation}) we can deduce
\[t_{p_j}(\varrho-\varphi)(x)-\frac{\log C_2}{A_{p_j}}-\left[\left(\frac{\varphi_{p_j}^{eq}-\varphi_{p_j}}{A_{p_j}}\right)-(\widetilde{\varphi^{eq}}-\varphi)\right](x)-\epsilon\leq \frac{\log P_{p_j}(x)}{2A_{p_j}}-\left(\frac{\varphi_{p_j}^{eq}-\varphi_{p_j}}{A_{p_j}}\right)(x).\]
Note $t_{p_j}(\varrho-\varphi)\to 0$ by definition and $[(\varphi_{p_j}^{eq}-\varphi_{p_j})/A_{p_j}-(\varphi^{eq}-\varphi)]\to 0$ in $L^1(X)$ by Lemma \ref{3.2.2}. So, it follows from (\ref{you need at end of this proof}) and the inequality above that

\[\limsup_j \left\|\frac{\log P_{p_j}}{2A_{p_j}}-\left(\frac{\varphi_{p_j}^{eq}-\varphi_{p_j}}{A_{p_j}}\right)\right\|_{L^1(X)}\leq \epsilon\int_X\frac{\omega^n}{n!}.\]
Letting $\epsilon\to 0$ tells us that

\[\frac{\log P_{p_j}}{2A_{p_j}}-\left(\frac{\varphi_{p_j}^{eq}-\varphi_{p_j}}{A_{p_j}}\right)\to 0\]
in $L^1(X),$ hence by (\ref{gamma-p curvature equation}) in the proof of Theorem \ref{T1.1} we have
\[\frac{\gamma_{p_j}-c_1(L_{p_j},h_{p_j}^{eq})}{A_{p_j}}\to 0\text{ as currents}.\]
From this and (\ref{T1.5-first-convergence}) it follows that 
\[\frac{c_1(L_{p_j},h_{p_j}^{eq})}{A_{p_j}}\to\theta +dd^c\varphi^{eq},\]
as desired.

\ii Next, we prove (\ref{alternate definition of varphi eq}). Since

\[\left(\frac{\log P_{p_j}}{2A_{p_j}}-(\varphi^{eq}-\varphi)\right),\left(\frac{\log P_{p_j}}{2A_{p_j}}-\frac{\varphi^{eq}_{p_j}-\varphi_{p_j}}{A_{p_j}}\right)\to 0\]
in $L^1(X)$, we also have 
\[\left(\frac{\varphi^{eq}_{p_j}-\varphi_{p_j}}{A_{p_j}}\right)\to (\varphi^{eq}-\varphi)\]
in $L^1(X)$. Then, from Lemma \ref{3.2.2} we deduce that $(\varphi^{eq}-\varphi)=(\widetilde{\varphi}^{eq}-\varphi)$ a.e., and $\widetilde{\varphi}^{eq}=\varphi^{eq}$ a.e. Since $\varphi^{eq}$ and $\widetilde{\varphi}^{eq}$ are qpsh and equal a.e., they are equal everywhere. Since $\varphi^{eq}$ is independent of the subsequence $\{p_j\},$ it follows from Lemma \ref{3.2.2} that every subsequence of $(\varphi_p^{eq}-\varphi_p)/A_p+\varphi$ has a subsequence converging to $\varphi^{eq}$ almost everywhere. Then $(\varphi_p^{eq}-\varphi_p)/A_p+\varphi\to \varphi^{eq}$ a.e., and (\ref{alternate definition of varphi eq}) holds, that is,
\[\varphi^{eq}=\left[\limsup_p\left(\frac{\varphi^{eq}_p-\varphi_p}{A_p}\right)\right]^*+\varphi.\]

\ii Lastly, we only proved the theorem for a particularly smooth $\theta\in \{T\}$, and we must show it holds for all. Let $\widetilde{\theta}\in \{T\}$ be smooth and choose $\widetilde{\varphi}\in C(X)$ such that $\widetilde{\theta}+dd^c \widetilde{\varphi}=T$. Set $\widetilde{\varphi^{eq}}=\sup\{\xi\in PSH(X,\widetilde{\theta})\divides \psi\leq \widetilde{\varphi}\}.$ Define $\sigma=\varphi-\widetilde{\varphi
}.$ We see that $\sigma \in C(X)$ and $dd^c\sigma= dd^c(\varphi-\widetilde{\varphi
})=\widetilde{\theta}-\theta,$
hence $\varphi^{eq}=\widetilde{\varphi^{eq}}+\sigma,$ and
\[\theta+dd^c\varphi^{eq}=(\widetilde{\theta}-dd^c\sigma)+dd^c(\widetilde{\varphi^{eq}}+\sigma)=\widetilde{\theta}+dd^c\widetilde{\varphi^{eq}}.\]
It follows that
\[\frac{\gamma_{p_j}}{A_{p_j}}\to\widetilde{\theta}+ dd^c\widetilde{\varphi^{eq}},\]
which shows the result is independent of our choice of $\theta$ as desired.
\eproof

\section{Applications}\label{S4}

Applying Theorem \ref{T1.2} to tensor products of powers of line bundles yields the following:

\begin{Corollary} Assume $F_j$, $1\leq j\leq k,$ are holomorphic line bundles on $X$ equipped with continuous metrics $h^{F_j}$ and singular metrics $g^{F_j}$ such that $c_1(F_j,g^{F_j})\geq 0$ and $c_1(F_1,g^{F_1})\geq \epsilon\omega,$ for some $\epsilon>0$. Let $\{m_{j,p}\}$ be a sequence of natural numbers with $m_{1,p}\to
 \infty, \text{ as }p\to\infty.$ Let $\gamma_p$, $p\geq 1,$ be the Fubini-Study currents associated with $H_{(2)}^0(X,L_p)$ where 
\[L_p=F_1^{m_{1,p}}\otimes\ldots\otimes F_k^{m_{k,p}}, \ h_p=(h^{F_1})^{m_{1,p}}\otimes\ldots\otimes (h^{F_k})^{m_{k,p}}.\]
Let $h_p^{eq}$ be the equilibrium metric of $h_p$. Then 
\[\frac{\gamma_p-c_1(L_p,h_p^{eq})}{\sum_{j=1}^k m_{j,p}}\to 0\]
weakly as currents.
\end{Corollary}

\emph{Proof.} Set
\smallskip
\[g_p=(g^{F_1})^{m_{1,p}}\otimes\ldots\otimes (g^{F_k})^{m_{k,p}}.\]

Let $x\in X$ and suppose $U,V\subset X$ are open polydisc neighborhoods of $x$ with $U\Subset V$. For $1\leq j\leq k$, let $\phi_p^{F_j}$ and $\rho_p^{F_j}$ denote the local weights of $h^{F_j}$ and $g^{F_j}$ on $V$, respectively. Then the local weights of $h_p$ and $g_p$ are given by
\smallskip
\[\phi_p:= \sum_{j=1}^k m_{j,p}\phi^{F_j}\text{ and }\rho_p:=\sum_{j=1}^k m_{j,p} \rho^{F_j}.\]

\ii Set
\smallskip
\[A_p:= \int_X c_1(L_p,g_p)\wedge \omega^{n-1}=\sum_{j=1}^k m_{j,p}\int_X c_1(F_j,g^{F_j})\wedge \omega^{n-1}.\]
Note there exists $c>0$, such that $cA_p \leq \sum_{j=1}^k m_{j,p}$
and
\[c_1(L_p,g_p)\geq m_{1,p}c_1(F_1,g^{F_j})\geq m_{1,p}\epsilon\omega.\]

By Theorem \ref{T1.2}, it suffices to show that $\{\phi_p/\sum_{j=1}^k m_{j,p}\}$ is equicontinuous and uniformly bounded in $L^1(U)$, and $\{\rho_p/\sum_{j=1}^k m_{j,p}\}$ is uniformly bounded in $L^1(U).$

\ii Observe that 
\[\frac{\phi_p}{\sum_{j=1}^k m_{j,p}}=\frac{\sum_{j=1}^k m_{j,p}\phi^{F_j}}{\sum_{j=1}^k m_{j,p}}=\sum_{j=1}^k \left(\frac{m_{j,p}}{\sum_{\ell=1}^k m_{\ell,p}}\right)\phi^{F_j},\]

so as $\{m_{j,p}/\sum_{\ell=1}^k m_{\ell,p}\}$ is bounded for all $j,$ and $\phi^{F_j}$ is continuous on $\overline{U},$ it follows that $\{\phi_p/\sum_{j=1}^k m_{j,p}\}$ is equicontinuous and uniformly bounded in $L^1(U).$ Similarly, 
\[\frac{\rho_p}{\sum_{j=1}^k m_{j,p}}= \sum_{j=1}^k \left(\frac{m_{j,p}}{\sum_{\ell=1}^k m_{\ell,p}}\right)\rho^{F_j},\]

hence $\{\rho_p/\sum_{j=1}^k m_{j,p}\}$ is uniformly bounded in $L^1(U).$ This completes the proof.
\eproof

\begin{Remark}

If $h^{F_j}=g^{F_j}$, $1\leq j\leq k$, then $h_p^{eq}=h_p$, and 

\[\frac{c_1(L_p,h_p)}{\sum_{\ell=1}^k m_{\ell,p}}=\sum_{j=1}^k\left(\frac{m_{j,p}}{\sum_{\ell=1}^k m_{\ell,p}}\right) c_1(F_j,h^{F_j}).\]

If additionally $m_{j,p}/\sum_{\ell=1}^k m_{\ell,p}$ converges for all $j$, then 
\smallskip
\[\frac{\gamma_p}{\sum_{\ell} m_{\ell,p}}\to \sum_{j=1}^k \left(\lim_p\frac{m_{j,p}}{\sum_{\ell} m_{\ell,p}}\right)c_1(F_j,h^{F_j}).\]

\end{Remark}

\ii The following results are corollaries of Theorem \ref{T1.4}:

 \begin{Corollary}\label{C1.4.1} Let $(X,\omega)$ be as in $(A)$ and for $p\geq 1$, let $L_p$ be a holomorphic line bundle on $X$ equipped with the metric $h_p$ with $C^2$ local weights $\phi^\alpha_p$ for all $\alpha.$ As well, let $A_p>0$ with $A_p\to\infty$. Suppose the following conditions hold. 

\begin{enumerate}[label=(\alph*)]

\item[(a)] There exists a continuous $(1,1)-$form $\Phi$ such that
\smallskip
\[\frac{c_1(L_p,h_p)}{A_p}\to \Phi\]

uniformly. 

\item[(b)] There exists $\varrho\in L^1(X)$ such that 
\smallskip
\[dd^c\varrho +\Phi\geq c\omega\]
\end{enumerate}

Then for any $\theta\in \{\Phi\}$, and $\varphi\in C(X)$ such that $dd^c\varphi=\Phi-\theta$, we have 
\[
\frac{\gamma_p}{A_p}\to \theta+dd^c\varphi^{eq}\text{ and }\frac{c_1(L_p,h_p^{eq})}{A_p}\to \theta+dd^c\varphi^{eq}
\]

weakly as currents, where $\varphi^{eq}=\sup\{\psi\in PSH(X,\theta)\divides \psi\leq\varphi\}.$

\end{Corollary}

\ii To prove the above Corollary, we will show that the given conditions imply the conditions stated in Lemma \ref{T1.4}. To do this, we first prove the following lemma.
\begin{Lemma}\label{3.2.1} Let $U$ be a coordinate polydisc. Assume $\{\phi_p:U\to \R\}$ is a family of $C^2$ functions bounded in $L^1$, and there exists a continuous $(1,1)-$form $\Phi$ such that 
\[dd^c\phi_p \to \Phi \text{ uniformly.}\]
Then, there exists $\phi\in C^1(U)$ and a subsequence $\{\phi_{p_k}\}$ such that
\[\phi_{p_k}\to \phi\text{ locally uniformly.}\]
\end{Lemma}

\emph{Proof.} Let $K\Subset U$ and $\epsilon>0$. Choose $r_0>0$ such that 
\[\widehat{K}=\{x\in U\divides \text{dist}(x,K)\leq r_0\}\Subset U.\]

Since $dd^c\phi_p \to \Phi$ uniformly it follows that there exists $M>0$ such that 
\[-M\omega\leq  dd^c \phi_p \leq  M\omega.\] 
Let $\zeta:U\to \R$ be a smooth real potential for $\omega$. As $\{\phi_p+M\zeta\}$ is a family of psh functions bounded in $L^1$, it follows from \cite[Theorem 3.2.12]{Hor94} that there exists $\widetilde{\phi}_U\in PSH(U)$ and a subsequence $\{\phi_{p_k}+M\zeta\}$ such that $\phi_{p_k}+M\zeta\to \widetilde{\phi}_U$ in $L^1_{loc}$. We define $\phi_U=\widetilde{\phi}_U-M\zeta$. We claim $\phi_U$ is continuous differentiable. 

\ii To see that $\phi_U$ is $C^1$ smooth we note that $\Delta \phi_U=\text{Tr}(\Phi)$ is continuous. If $K\Subset U$ is a ball, then $\Delta \phi_U$ is bounded on $K$, so by \cite[Lemma 4.1]{GT98} there exists $\psi_K\in C^1(K)$ with $\psi_K=\phi_U$ a.e. on $K.$ Then, $\psi_K+M\zeta=\phi_U+M\zeta$ a.e. and both sides are psh, so they're equal everywhere. It follows that $\phi_U\in C^1(K)$. As $K\subseteq U$ was an arbitrary ball we have $\phi_U\in C^1(U)$.

\ii Next, as $dd^c\phi_{p_k} \to\Phi$ uniformly, it follows that for $k$ large enough, we have
\[-\omega\leq dd^c (\phi_{p_k}-\phi_U)\leq \omega.\]
From this and subaveraging, we find
\[(\phi_{p_k}-\phi_U+\zeta)(x)< \frac{n!}{\pi^n r^{2n}}\int_{B(x,r)}(\phi_{p_k}-\phi_U+\zeta)d\lambda,\]
for all $x\in K$ and $0<r<r_0$, where $d\lambda$ is the Lebesgue measure. It follows that
\[(\phi_{p_k}-\phi_U)(x)\leq\frac{n!}{\pi^n r^{2n}}\int_{B(x,r)}(\phi_{p_k}-\phi_U)d\lambda+\frac{n!}{\pi^n r^{2n}}\int_{B(x,r)}\left(\zeta-\zeta(x)\right)d\lambda.\]
By the uniform continuity of $\zeta$ on $\widehat{K}$, it follows that we can choose $r$ small enough so that
\[(\phi_{p_k}-\phi_U)(x)\leq\frac{n!}{\pi^n r^{2n}}\int_{B(x,r)}(\phi_{p_k}-\phi_U)d\lambda+\epsilon\]
for all $x\in K$. By applying similar arguments to $(\phi_U-\phi_{p_k})$, it follows that by shrinking $r$, if necessary, we may assume
\[(\phi_U-\phi_{p_k})(x)\leq \frac{n!}{\pi^n r^{2n}}\int_{B(x,r)}(\phi_U-\phi_{p_k})d\lambda+\epsilon.\]
Combining these inequalities gives
\[\frac{n!}{\pi^n r^{2n}}\int_{B(x,r)}(\phi_{p_k}-\phi_U)d\lambda-\epsilon\leq (\phi_{p_k}-\phi_U)(x)\leq\frac{n!}{\pi^n r^{2n}}\int_{B(x,r)}(\phi_{p_k}-\phi_U)d\lambda+\epsilon,\]
so 
\[|\phi_{p_k}-\phi_U|(x)\leq\frac{n!}{\pi^n r^{2n}}\int_{B(x,r)}|\phi_{p_k}-\phi_U|d\lambda+\epsilon\leq \frac{n!}{\pi^n r^{2n}}\int_{\widehat{K}}|\phi_{p_k}{-\phi_U}|d\lambda+\epsilon.\]
As $\phi_{p_k}\to \phi_U$ in $L^1(\widehat{K})$, it follows that if $k$ is large enough, then
\[|\phi_{p_k}-\phi_U|(x)\leq \frac{n!\hspace{.3mm}\epsilon}{\pi^n r^{2n}}+\epsilon.\]
As $x$ was arbitrary, we see $\phi_{p_k}\rightarrow \phi_U$ uniformly on $\widehat{K},$ therefore $\{\phi_p\}$ has a subsequence converging uniformly to $\phi_U$ on $K$.

\ii To see that $\{\phi_p\}$ has a subsequence converging locally uniformly on $U,$ we let $\{K_n\}$ be an exhaustion of $U$ by compact sets, and then apply a diagonalization argument to construct a sequence converging uniformly on each $K_n$, and hence, locally uniformly on $K.$
\eproof

\textbf{Proof of Corollary \ref{C1.4.1}:} As $dd^c\phi_p^\alpha/A_p\to \Phi$ uniformly on $V_\alpha$ it follows that there exists $M>0$ such that 
\[-M\omega\leq dd^c\phi^\alpha_p/A_p\leq M\omega.\]
Hence, there exists a smooth function $\zeta$ such that $\{\phi^\alpha_p/A_p+\zeta\}$ is plurisubharmonic for all $p,\alpha$. By \cite[Proposition A.16]{DS98}, it follows that there exists a $c>0$ and $\widehat{\phi^\alpha_p}\in PSH(U_\alpha)$ for all $p,\alpha$ such that $dd^c\widehat{\phi^\alpha_p}=dd^c(\phi_p^\alpha+A_p\zeta)$ and 
\[\|\widehat{\phi_p^\alpha}/A_p \|_{L^1(U_\alpha)}\leq c\int_{V_\alpha}dd^c(\phi_p^\alpha/A_p+\zeta)\wedge\omega^{n-1}.\] 
Define $\widetilde{\phi^\alpha_p}/A_p=\widehat{\phi^\alpha_p}/A_p-\zeta$. From above it is clear that $\{\widetilde{\phi^\alpha_p}/A_p\}$ is bounded in $L^1(U_\alpha)$. As $dd^c \widetilde{\phi^\alpha_p}=dd^c \phi_p^\alpha$, it follows that $\widetilde{\phi^\alpha_p}=\phi_p^\alpha+k_p^\alpha$ almost everywhere for some pluriharmonic function $k_p^\alpha\in C^\infty (U_\alpha)$. Hence, $\widetilde{\phi^\alpha_p}$ is $C^2$, and $\widetilde{\phi^\alpha_p}=\phi_p^\alpha+k_p^\alpha$ everywhere. 

\ii Since $k_p^\alpha$ is pluriharmonic, it is the real part of a holomorphic function $f_p^\alpha$. Let $\widetilde{e_p^\alpha}=e_p^\alpha e^{-f_p^\alpha}$. We find 
\[h_p(\widetilde{e_p^\alpha},\widetilde{e_p^\alpha})=e^{-f_p^\alpha}\overline{e^{-f_p^\alpha}}h_p(e_p^\alpha,e_p^\alpha)=e^{-2\text{Re}f_p^\alpha}e^{-2\phi_p^\alpha}=e^{-2(\phi_p^\alpha+k_p^\alpha)}=e^{-2\widetilde{\phi_p^\alpha}},\]
so $\widetilde{\phi_p^\alpha}$ are $C^2$ local weights of $h_p$ with respect to the local frames $\widetilde{e_p^\alpha}$, and the family $\{\widetilde{\phi_p^\alpha}/A_p\}$ is bounded in $L^1(U_\alpha).$

\ii Suppose $\{p_j\}$ is a subsequence. By condition (a) and Lemma \ref{3.2.1}, it follows that for all $\alpha$, there exists $\phi^\alpha\in C^1(U_\alpha)$ and a subsequence $\{p_{j_k}\}$ such that $\widetilde{\phi_{p_{j_k}}^\alpha}/A_{p_{j_k}}\to \phi^\alpha$ locally uniformly on $U_\alpha$. We will start by showing that the claim holds for this subsequence.

\ii We shall apply Theorem \ref{T1.4} to the subsequence $\{p_{j_k}\}$ with open cover $\{U_\alpha\}$ and local frames $\widetilde{e_p^\alpha}$. Given $\widetilde{\phi_{p_{j_k}}^\alpha}/A_{p_{j_k}}\to \phi^\alpha$, conditions $(A)$ and $(B')$ are satisfied for our subsequence. To see that (b) of Theorem \ref{T1.4} holds, we first note that $dd^c(\varrho+\phi^\alpha)=dd^c\varrho+\Phi|_{V_\alpha}=dd^c\varrho+ T|_{V_\alpha}.$ Next, we define $\varrho'$ by $\varrho'=\varrho+\varphi$. Then,
\[dd^c(\varrho+\phi^\alpha)=dd^c\varrho+ T|_{V_\alpha}=dd^c\varrho+\theta+dd^c\varphi=dd^c\varrho'+\theta,\]
so $dd^c\varrho'+\theta\geq c\omega$, and (b) is satisfied. It remains to show that condition $(a)$ of Theorem \ref{T1.4} holds for the subsequence. 

\ii As $c_1(L_p,h_p)/A_p$ is uniformly convergent to $\Phi=dd^c\phi^\alpha$, it follows that there exists $\delta_p\geq0 $ with $\delta_p\to 0$ such that
\[dd^c \phi^\alpha-\delta_p\omega \leq \frac{c_1(L_p,h_p)}{A_p},\]
which gives us (a) of Theorem \ref{T1.4}.

\ii It follows from Theorem \ref{T1.4} that for any $\theta\in \{\Phi\}$ and $\varphi\in C(X)$ such that $dd^c\varphi=\Phi-\theta$, we have
\[\frac{\gamma_{p_{j_k}}}{A_{p_{j_k}}}\to\theta+dd^c \varphi^{eq}\text{ and } \frac{c_1(L_{p_{j_k}},h_{p_{j_k}}^{eq})}{A_{p_{j_k}}}\to\theta+dd^c \varphi^{eq},\]
where 
\[
    \varphi^{eq}=\sup\{\psi\in PSH(X,\theta) \divides  \psi\leq \varphi \}.
\]
Since this convergence holds for an arbitrary subsequence, it holds for the entire sequence.
\eproof
\begin{Remark} In the above scenario, we do get for all $\alpha$ that $\phi^\alpha_p/A_p$ has a subsequence converging uniformly to some $\phi^\alpha\in C^1(U_\alpha)$. Unlike the $\varphi^{eq}$ in our proof above, the function $\phi^\alpha$ is dependent on the subsequence. However, in the case where $L_p=L^p$, the function $\phi^\alpha$ is unique up to constant addition. 
\end{Remark}
\ii Before stating our next corollary, we let $\{V_\alpha\}$ be an open cover of $X$ as in Section \ref{S2.2}, and define $\phi_p^\alpha,\rho_p^\alpha$ to be the local weights of $h_p$ and $g_p$ on $V_\alpha.$
\begin{Corollary}\label{C1.4.2}Let $(X,\omega)$ and $L_p$ be as in $(A)$ and $(B)$. Assume the local weights $\phi_p^\alpha$ of $h_p$ are $C^2$ and $\{\phi_p^\alpha/A_p\}$ is bounded in $L^1(V_\alpha)$. Suppose the following conditions hold. 

\begin{enumerate}[label=(\alph*)]

\item[(a)] There exists a continuous $(1,1)-$form $\Phi$ such that
\smallskip
\[\frac{c_1(L_p,h_p)}{A_p}\to \Phi\]

uniformly. 

\item[(b)] For all $\alpha$ the family $\{\rho_p^\alpha/A_p\}$ is uniformly bounded in $L^1(V_\alpha)$, and there exists $C\geq 0$ such that $A_p/a_p\leq C.$

\end{enumerate}
Then for any real smooth closed $\theta\in \{\Phi\}$, and $\varphi\in C(X)$ such that $dd^c\varphi=\Phi-\theta$, we have 
\[
\frac{\gamma_p}{A_p}\to \theta+dd^c\varphi^{eq}\text{ and }\frac{c_1(L_p,h_p^{eq})}{A_p}\to \theta+dd^c\varphi^{eq}
\]

weakly as currents, where $\varphi^{eq}=\sup\{\psi\in PSH(X,\theta)\divides \psi\leq\varphi\}.$
\end{Corollary}

\textbf{Proof:} Note that, since $(B)$ is assumed, our setup is stronger than the one in Corollary \ref{C1.4.1}, so (a) implies condition (a) of Corollary \ref{C1.4.1}. It suffices to show that condition (b) of Corollary \ref{C1.4.1} holds, i.e., there exists a $\varrho\in L^1(X)$ such that $\Phi+dd^c\varrho>c\omega$ for some $c>0.$

\ii By Lemma \ref{3.2.1}, there exists $\phi^\alpha\in C^1(V_\alpha)$ and a subsequence $\{\phi_{p_j}^\alpha\}$ such that $\phi_{p_j}^\alpha/A_p\to \phi^\alpha$ locally uniformly. Let $\zeta$ be a real smooth potential for $\omega$ on $V_\alpha$. As $A_p/a_p$ is bounded, it follows that for some $c>0$ we have
\smallskip
\[dd^c(\rho_p^\alpha/A_p-c\zeta)=c_1(L_p,g_p)/A_p-c\omega\geq (a_p/A_p)\omega-c\omega\geq 0,\]
 
 so $(\rho_p^\alpha/A_p-c\zeta)$ is psh for all $\alpha.$ Since $(\rho_p^\alpha/A_p-c\zeta)$ is uniformly bounded in $L^1_{loc}(V_\alpha)$, it follows from \cite[Theorem 3.2.12]{Hor94} that $(\rho_{p_j}^\alpha/A_{p_j}-c\zeta)$ has a subsequence converging to a function $\widetilde{\rho}^\alpha\in PSH(U_\alpha)$ in $L^1(U_\alpha)$. By refining our original subsequence, we may assume WLOG that the convergence holds for all $\alpha$ simultaneously.

\ii We define $\rho^\alpha=\widetilde{\rho}^\alpha+c\zeta,$ and set $\varrho^\alpha=\rho^\alpha-\phi^\alpha.$ We compute
\smallskip
\[
\Phi+dd^c\varrho^\alpha=\Phi+dd^c\rho^\alpha-dd^c \phi^\alpha=dd^c(\widetilde{\rho}^\alpha)+dd^c c\zeta\geq c\omega. 
\]

Next, notice that
\smallskip
\[\frac{\rho^\alpha_{p_j}-\phi^\alpha_{p_j}}{A_{p_j}}\to \rho^\alpha-\phi^\alpha=\varrho^\alpha\]

in $L^1(U_\alpha)$ for all $\alpha$. The left-most function glues to a global function, hence $\varrho^\alpha=\varrho^\beta$ almost everywhere for all $\alpha,\beta$. Since $\varrho^\alpha$ is qpsh for all $\alpha$, the equality holds everywhere. Then, the function $\varrho:X\to[-\infty,\infty)$ defined by $\varrho|_{V_\alpha}=\varrho^\alpha$ is well defined and
satisfies
\[\Phi+dd^c \varrho\geq c\omega.\]
Clearly, $\varrho\in L^1(x)$, so condition (b) of Corollary \ref{C1.4.1} holds.
\eproof

\ii Before stating the next theorem we recall than in $(B')$, the functions $\phi^\alpha$ were the local uniform limits of the scaled local weights $\phi^\alpha_p/A_p$. 
\begin{Theorem}\label{T1.5} Let $(X,\omega)$ be as in (A) and $L$ be a big line bundle on $X$. Suppose $(L_p,h_p)$ are as in $(B')$, where $L_p=L^p$, $e_p^\alpha=(e_1^\alpha)^{\otimes p}$, and $A_p=p$. Then, there is a metric $h$ on $L$ with local weights $\phi^\alpha$ on $V_\alpha,$ and

\begin{equation}\label{T1.5 convergence}
\frac{\gamma_p}{p}\to c_1(L,h^{eq})\text{ and }\ \frac{c_1(L^p,h_p^{eq})}{p}\to c_1(L,h^{eq})\text{ weakly as currents},
\end{equation}
where $h^{eq}$ is the equilibrium metric determined by $h$, as defined in Section \ref{S2.2}.
\end{Theorem}

\ii The convergence (\ref{T1.5 convergence}) is a special case of (\ref{T1.4-convergence}), and one can formulate a statement similar to (\ref{T1.5 convergence}) whenever the $\phi^\alpha$ are local weights of a metric. This isn't the only similarity to Theorem \ref{T1.4}. In fact, the only condition from  Theorem \ref{T1.4} which we can't show holds here is the assumption that $c_1(L_p,h_p)/A_p\geq dd^c\phi^\alpha-\delta_p\omega$. This property isn't needed here, as $g$ handles the positivity requirement. Despite the similarities, Theorem \ref{T1.5} is proved using Theorem \ref{T1.2} and not Theorem \ref{T1.4}.

\ii The following lemma will be needed in our proof of Theorem \ref{T1.5}. 

\begin{Lemma}\label{Envelope-Lemma} If $\theta$ is a real smooth closed $(1,1)-$form on $X$ with $PSH(X,\theta)\not=\emptyset$, and $\sigma_p\in C(X)$ converges uniformly to $\sigma\in C(X)$, then $\tau_p,\tau\in PSH(X)$ with $\tau_p\rightarrow\tau$ uniformly, where
\[\tau_p:=\sup\{\psi\in PSH(X,\theta)\divides  \psi\leq \sigma_p\}\text{ and } \tau:=\sup\{\psi\in PSH(X,\theta)\divides \psi\leq \sigma\}.\]
\end{Lemma}
\emph{Proof.} Define
\[F_p=\{\psi\in PSH(X,\theta)\divides  \psi\leq \sigma_p\}\text{ and }F=\{\psi\in PSH(X,\theta)\divides \psi\leq \sigma\}.\]

Let $\psi\in F_p$. We see 
\[\psi-\|\sigma_p-\sigma\|_{L^\infty}\leq \sigma\text{ and }\psi-\|\sigma_p-\sigma\|_{L^\infty}\in PSH(X,\theta),\]

so $\psi-\|\sigma_p-\sigma\|_{L^\infty}\in F.$ Taking the supremum over all such $\psi$ we obtain $\tau_p-\|\sigma_p-\sigma\|_{L^\infty}\leq \tau.$ Similarly, $\tau-\|\sigma-\sigma_p\|_{L^\infty}\leq \tau_p,$ thus $\|\tau-\tau_p\|_{L^\infty}\leq\|\sigma_p-\sigma\|_{L^\infty},$ and so $\tau_p\rightarrow\tau$ uniformly as desired.
\eproof

\ii For the following proof, we fix a smooth metric $h_0$ on $L$ with curvature $\theta=c_1(L,h_0)$. Recall definition (\ref{def of T current}), that is, $T$ is the current satisfying $T|_{V_\alpha}=dd^c\phi^\alpha$. By definition, the global weight of $h^{eq}$ is given by $\varphi^{eq}=\sup\{\psi\in PSH(X,\theta)\divides \psi\leq\varphi\}.$

\smallskip
\textbf{Proof of Theorem \ref{T1.5}.} As in the proof of Proposition \ref{Prop}, there exist pluriharmonic functions $\psi_p^{\alpha\beta}\in C^\infty(V_\alpha\cap V_\beta)$ for all $\alpha,\beta$, and $p$, which are defined by $\psi^{\alpha\beta}_p=\log|k^{\alpha\beta}_p|,$ where $k^{\alpha\beta}_p$ denotes the transition function of $L_p$. These have the property that
\smallskip
\[\phi_p^\alpha=\phi_p^\beta+\psi_p^{\alpha\beta}.\]

Since $e_p^\alpha=(e^\alpha)^{\otimes p}$, it follows that  $\psi_p^{\beta\alpha}=p\psi_1^{\beta\alpha}$, thus 
\smallskip
\[\phi^\alpha=\lim_p \left(\frac{\phi_p^\alpha}{p}\right)=\lim_p\left(\frac{\phi_p^\beta+p\psi_1^{\beta\alpha}}{p}\right)=\phi^\beta+\psi_1^{\beta\alpha},\]

hence $h$ is well defined, and $\phi_p^\alpha/p-\phi^\alpha$ glue to a global potential for $c_1(L_p,h_p)/p-T.$

\ii Let $g$ be a metric on $L$ with strictly positive curvature current and local weights $\rho^\alpha$. The family of bounded local weights of $g^p$ is $\{p\cdot \rho^\alpha/A_P\}=\{\rho^\alpha\},$ which is clearly bounded in $L^1(V_\alpha)$. Since $\phi_p^\alpha/p\rightarrow \phi^\alpha$ locally uniformly on $V_\alpha$, it follows that $\{\phi_p^\alpha/p\}$ is equicontinuous and uniformly bounded on $U_\alpha$. This shows the hypothesis of Theorem \ref{T1.2} hold with $g_p=g^p$, hence
 \smallskip
 \[\frac{\gamma_p-c_1(L_p,h_p^{eq})}{p}\to 0.\]
Since $c_1(L,h^{eq})=\theta+dd^c \varphi^{eq}$, we'll be done if we show
\[\frac{c_1(L_p,h_p^{eq})}{p}\to \theta+dd^c \varphi^{eq}.\]
As
\smallskip
\[\frac{c_1(L_p,h_p^{eq})}{p}=\theta+\frac{ dd^c \varphi_p^{eq}}{p},\]

our desired convergence will hold provided
\smallskip
\[\frac{\varphi_p^{eq}}{p}\rightarrow \varphi^{eq}\text{ locally uniformly}.\]
We show this now.

\ii  Let $\varphi$ and $\varphi_p$ denote the global weights of $h$ and $h_p$ with respect to $h_0$ and $h_0^p$, respectively. If $\psi^\alpha$ denotes the local weight of $h_0$ on $V_\alpha$, then 
\[\frac{\varphi_p}{p}=\frac{\phi_p^\alpha-p\psi^\alpha}{p}=\frac{\phi_p^\alpha}{p}-\psi^\alpha.\]
Then ${\varphi_p}/{p}\to \phi^\alpha-\psi^\alpha=\varphi$ uniformly. By definition we have
\[\frac{\varphi_p^{eq}}{p}=\sup\{\psi\in PSH(X,\theta)\in \psi \leq \varphi_p/p\}\text{ and }{\varphi^{eq}}=\sup\{\psi\in PSH(X,\theta)\in \psi \leq \varphi\},\]
hence by Lemma \ref{Envelope-Lemma} it follows that $\varphi_p^{eq}/p\to \varphi^{eq}$ locally uniformly, which completes the proof
\eproof

\end{document}